\documentclass{amsart}
\usepackage[utf8]{inputenc}

\usepackage{amssymb}
\usepackage{amsfonts}
\usepackage{amsmath}
\usepackage{amstext}
\usepackage{amsthm}
\usepackage{latexsym}

\usepackage[utf8]{inputenc}

\usepackage{mathrsfs}  
\usepackage[mathscr]{eucal}

\usepackage{amscd}
\usepackage{eucal}
\usepackage{amsxtra}
\usepackage{amsbsy}
\usepackage{graphicx}
\usepackage{latexsym}
\usepackage{pb-diagram}
\usepackage{amsopn}
\usepackage{delarray}
\usepackage{psfrag}
\usepackage{verbatim}
\usepackage{bussproofs}
\usepackage[T1]{fontenc}
\usepackage{color,cancel}

\usepackage{appendix}
\usepackage{color}
\usepackage[dvipsnames]{xcolor}
\usepackage{epigraph}
\usepackage{lipsum}
\usepackage{hyperref} 
\hypersetup{pdfpagemode={UseOutlines},
bookmarksopen=true,bookmarksopenlevel=0,hypertexnames=false,colorlinks=true,citecolor=black,linkcolor=black,urlcolor=black,pdfstartview={FitV},unicode,breaklinks=true}
\usepackage{graphicx}
\usepackage{enumitem}

\newtheorem{mydef}{Definition}%
\newtheorem{theorem}[mydef]{Theorem}
\newtheorem{lemma}[mydef]{Lemma}
\newtheorem{notation}[mydef]{Notation}
\newtheorem{proposition}[mydef]{Proposition}
\newtheorem{corollary}[mydef]{Corollary}

\newtheorem{example}[mydef]{Example}

\newtheorem{question}[mydef]{Question}

\theoremstyle{remark}
\newtheorem{claim}[mydef]{Claim}

\def\cL{{\mathcal{L}}}
\def\cI{{\mathcal{I}}}
\def\cJ{{\mathcal{J}}}
\def\cF{{\mathcal{F}}}
\def\cP{{\mathcal{P}}}

\def\cS{{\mathcal{S}}}

\def\cH{{\mathcal{H}}}

\def\cM{{\mathcal{M}}}

  \def\RR{{\mathbb{R}}}

\def\cA{{\mathcal{A}}}
\def\cB{{\mathcal{B}}}
\def\posit{{\mathcal{I}^+(\mathcal{A})}}
\def\dobleposit{{\mathcal{I}^{++}(\mathcal{A})}}

\def\s{{\mathfrak{s}}}
\def\a{{\mathfrak{a}}}
\def\b{{\mathfrak{b}}}
\def\c{{\mathfrak{c}}}
\def\p{{\mathfrak{p}}}
\def\d{{\mathfrak{d}}}

\def\r{{\mathfrak{r}}}

\newcommand{\conc}{^\smallfrown}
\newcommand{\N}{\omega}
\newcommand{\w}{\omega}
\newcommand{\cont}{\mathfrak{c}}

\newcommand{\rest}{\upharpoonright}

\newcommand{\Ramsey}{sequentially compact }
\newcommand{\RRamsey}{sequentially compact}

\def\cL{{\mathcal{L}}}

\title{High dimensional sequential compactness}

\author{C\'esar Corral}
\address{York University, 4700 Keele St, Toronto, ON, Canada, M3J 1P3.}
\email{cicorral@yorku.ca}

\author{Osvaldo Guzm\'an}
\address{Centro de Ciencias  Matem\'aticas, Universidad Nacional Aut\'onoma de M\'exico,
\' Campus Morelia, 58089, Morelia, Michoac\'an, M\'exico.}
\email{oguzman@matmor.unam.mx}
\thanks{The research of the second author was supported by PAPIIT grant IA102222.}

\author{Carlos L\'opez-Callejas}
\address{Centro de Ciencias  Matem\'aticas, Universidad Nacional Aut\'onoma de M\'exico,
\' Campus Morelia, 58089, Morelia, Michoac\'an, M\'exico.}
\email{carloscallejas@matmor.unam.mx}
\thanks{The research of the third author was supported by PAPIIT grant IN101323 and CONACyT grant A1-S-16164.}

\keywords{$n$-\Ramsey space, sequentially compact, almost disjoint family, splitting number, bounding number.}

\subjclass[2020]{54A20, 03E02, 03E17, 54D80, 54D30}

\begin{document}

\maketitle

\begin{abstract}
    We give examples of $n$-\Ramsey spaces that are not $(n+1)$-\Ramsey under several assumptions. We improve results from Kubis and Szeptycki by building such examples from $\mathfrak{b=c}$ and $\diamondsuit(\mathfrak{b})+\mathfrak{d}=\omega_1$. We also introduce a new splitting-like cardinal invariant and then show that the same holds under $\mathfrak{s=b}$. 
\end{abstract}

\section{Introduction}

A widely studied topological property is that of sequential compactness. Recall that a space $X$ is \emph{sequentially compact} if every sequence $\{x_n:n\in\w\}\subseteq X$ has a convergent subsequence. 

An important fact about this property is that it characterizes compactness in metric spaces, that is, a metric space is sequentially compact if and only if it is compact.
Nevertheless, this fact is not true in general, since $\omega_1$ is a sequentially compact non compact space and $\beta\w$, the Stone-C\v{e}ch compactification of a countable discrete space, is compact but fails to be sequentially compact (see \cite{willard} ).

In \cite{kubis2021topological}, Kubi\'s and Szeptycki introduced higher dimensional versions of sequential compactness and showed that these properties also hold in metric spaces.
The $2$-dimensional case had been previously introduced and proved to hold in compact metric spaces in \cite{ramseymetric} by Boja\'{n}czyk, Kopczy\'{n}ski and Toru\'{n}czyk. There the authors used this to obtain idempotents in compact semigroups.

We follow the convention that $[B]^n$ denotes the family of all subsets of $B$ of cardinality $n$. Analogously, $[B]^\w$ denotes the family of all of its countable infinite subsets and $[B]^{<\w}=\bigcup_{n\in\w}[B]^n$. Recall that if $X$ is a topological space, $p\in X$, $B$ is a countable set and $f:B\rightarrow X$, then we say that $f$ \emph{converges to} $p$ if for all $U\subseteq X$ open such that $p\in U$, there exists a finite set $F\subseteq B$ such that $f''(B\setminus F)\subseteq U$. It is natural to generalize this definition to higher dimensions.
Following \cite{kubis2021topological}, if $n\geq 1$ and $f:[B]^n\rightarrow X$, then $f$ \emph{converges} to $p$ if for every neighborhood $U$ of $p$ there is a finite set $F\subseteq B$ such that $f''[B\setminus F]^n\subseteq U$ and in this case we say that the function is $f$ \emph{convergent}\footnote{Remind that $f''A=\{f(x):x\in A\}$.}.

\begin{mydef}\cite{kubis2021topological}
Let $X$ be a topological space and $n\geq 1$. We say that $X$ is an $n$-\Ramsey space if for every function $f:[\w]^n\rightarrow X$ there exists $B\in[\w]^{\w}$ such that $f\restriction [B]^n$ is convergent. We say that $f\restriction [B]^n$ is a convergent subsequence of $f$.
\end{mydef}

Observe that being $1$-\Ramsey is equivalent to being sequentially compact. It is also not hard to see that Ramsey's Theorem\footnote{Ramsey's theorem states that if  $n,m\geq 1$ and $f:[\mathbb{N}]^n\rightarrow m$, then there exists $Y\subseteq \mathbb{N}$ infinite such that $|f''[Y]^n|=1$.} is equivalent to the assertion that every finite space $X$ is $n$-\Ramsey  for all $n\geq 1$. 

A similar notion was studied by Helmut Knaust in \cite{angelicspaces} and \cite{arrayconvergence}, where the author defined Ramsey-uniformly convergence and making use of it, he defined a space to have the \emph{Ramsey property} if any double sequence $\{x(i,j):i,j\in\omega\}$ such that $\lim_{i\to\infty}\lim_{j\to\infty}=x$ also converges Ramsey-uniformly to $x$. The notion of Ramsey-uniformly convergence coincides with the notion of convergence introduced in \cite{kubis2021topological} in the 2-dimensional case. It was shown by Knaust that any pointwise compact subset of the Baire-$1$ functions on a Polish space has the Ramsey property \cite{arrayconvergence}, that any angelic space so that $\chi(x,A)<\p$ for all countable relatively compact subsets $A$ in $X$ and all $x\in X$ has the Ramsey property and that so do all spaces of continuous functions $C(Y)$ with the topology of pointwise convergence where $Y$ is a quasi Suslin space \cite{angelicspaces}. This result in the particular case of first countable spaces (and then for the real numbers with their usual topology) was previously showed by Rosenthal (see \cite{someremarks}).

In \cite{kubis2021topological}, the term $n$-Ramsey was used instead of $n$-sequentially compact, however, it seems that the latter is more accurate in this context as the term $n$-Ramsey conflicts with other previously introduced notions and does not stress the compactness nature of the property.
Note that the Ramsey property and being $2$-Ramsey (what we are calling now 2-sequentially compact) are similar but different. While the Ramsey property ensures a Ramsey-like convergence of a double sequence that already comes with a limit point, 2-sequential compactness also claims the existence of such a point $x\in X$ for every function $f:[\w]^2\to X$. For this reason and to avoid confusion with Knaust's property and Ramsey spaces introduced by Todor\v{c}evi\'c \cite{Todorcevic+2010}, we have decided to adopt the more descriptive name of $n$-sequential compactness

In \cite{kubis2021topological}, Kubi\'s and Szeptycki proved the following result:

\begin{theorem}\cite{kubis2021topological}
    Every compact metric space is $n$-\Ramsey for every $n\geq 1$.
\end{theorem}

Thus we can not differentiate the class of $n$-\Ramsey spaces from the class of $m$-\Ramsey spaces with $n\neq m$ in the realm of metric spaces. So, the interesting question is whether we can distinguish these properties in the general category of topological spaces.

It is easy to see that every $(n+1)$-\Ramsey space is also $n$-\Ramsey (Proposition 2.2 in \cite{kubis2021topological}) and Kubi\'s and Szeptycki proved the following:

\begin{theorem}\label{teo3}
\cite{kubis2021topological}
    \begin{enumerate}
        \item\label{item1} There is a compact sequentially compact space (that is, $1$-\RRamsey) that is not $2$-\RRamsey.
        
        \item\label{item2} Under \emph{\textsf{CH}}, for every $n\in\omega$, there exists a $n$-\Ramsey space that is not $(n+1)$-\RRamsey. 
    \end{enumerate}
\end{theorem}



In this work we will improve item (\ref{item2}) of Theorem \ref{teo3} by constructing $n$-\Ramsey spaces that are not $(n+1)$-\Ramsey under several cardinal invariant assumptions. We also show that there is Fr\'echet example for item (\ref{item1}) answering a question posed by the authors in \cite{kubis2021topological}. This example was sketched in that paper under the assumption of a completely separable mad family, which follows under, e.g., $\mathfrak{c}\leq\aleph_\w$. We point out that the construction does not need a maximal almost disjoint family and since the existence of completely separable almost disjoint families follows from \textsf{ZFC}, an easy modification of the result in \cite{kubis2021topological} can be carried out without extra assumptions. Although we will be quoting theorems from \cite{kubis2021topological}, it is not needed that the reader is familiar with it, since we will be recalling the needed notions and results.

The theory of $n$-\Ramsey spaces is just in its beginnings and we expect it to be richer and more interesting than the one of sequentially compact spaces. This is our humble contribution to this new theory.

\section{Preliminaries}\label{preliminaries}

All spaces considered here are Hausdorff. Denote by $\omega$ the first infinite ordinal, by $\omega_1$ the first uncountable ordinal and by $\c$ the cardinality of $\RR$.
Given two infinite sets $A,B\subseteq\w$, we write $A\subseteq^* B$ if $A\setminus B$ is finite. Two sets $A,B$ are \emph{almost disjoint} if $A\cap B$ is finite. A family $\cA\subseteq[\w]^\w$ is \emph{almost disjoint (ad)} if its elements are pairwise almost disjoint and it is \emph{maximal almost disjoint (mad)} if it is ad and is not properly contained in any other almost disjoint family. 

We say that a set $S\in[\w]^\w$ \emph{splits} $A\subseteq\w$ if $|S\cap A|=\w=|A\setminus S|$. A family $\cS\subseteq[\w]^\w$ is a \emph{splitting family} if for every $X\in\omega$ there exists $S\in\cS$ such that $S$ splits $X$. 

The set of functions from $X$ to $Y$ will be denoted by $X^Y$. For two functions $f,g\in\omega^\w$, we write $f\leq^* g$ if $\{n\in\omega:f(n)>g(n)\}$ is finite. A family $\cB\subseteq\w^\w$ is \emph{unbounded} if for every $f\in\w^\w$ there exists $b\in\cB$ such that $b\nleq^*f$ and $\cB$ is \emph{dominating} if for every $f\in\w^\w$ there exists $d\in \cB$ such that $f\leq^*d$. We also say that $\cB$ is a scale if $\cB$ is a dominating family and it is well ordered by $\leq^*$.

We define below the cardinal invariants used in this work for the ease of the reader,
\begin{itemize}
    \item $\a$ is the minimum size of a mad family,
    \item $\s$ is the minimum size of a splitting family,
    \item $\b$ is the minimum size of an unbounded family,
    \item $\d$ is the minimum size of a dominating family.
\end{itemize}

Each of these cardinals is uncountable and smaller or equal to $\c$.
It is also well known that $\b\leq\d$, $\b\leq\a$ and $\s\leq\d$ and there are no other provable relations among them in $\textsf{ZFC}$.
It is worth noting that scales exist if and only if $\b=\d$. A nice survey about cardinal characteristics and their properties is \cite{Blass2010}.

A family $\cI\subseteq\cP(\omega)$ is an \emph{ideal} if it is closed under finite unions and subsets, i.e., if $A,B\in\cI$ then $A\cup B\in\cI$ and if $A\subseteq B$ and $B\in \cI$ then $A\in\cI$. An ideal $\cI$ is \emph{tall} if for every $X\in[\w]^\w$ there exists $I\in\cI$ such that $|X\cap I|=\w$. Given a tall ideal $\cI$, the cardinal $\textsf{cov}^*(\cI)$ is defined as follows:
\[
\textsf{cov}^*(\cI)=\min\{|\cH|: \cH\subseteq\cI\  (\forall A\in [\w]^\w\exists B\in\cH (|A\cap B|=\w))\}.
\]

We can also define $\textsf{cof}(\cI)$, the cofinality of an ideal $\cI$, as the minimum size of a subfamily of $\cI$ that generates it, that is:
\[\textsf{cof}(\cI)=\min\{|\cJ|:\cJ\subseteq\cI\land\forall\ I\in\cI\ \exists J\in\cJ(I\subseteq J)\}.\]

For more on ideals on countable sets and their cardinal invariants the reader may consult \cite{hrusakcombinatoricsoffiltersandideals}, \cite{cardinalinvariantsofanalytic}, \cite{pairsplitting}, \cite{densityzeroideal} and \cite{twoinequalities}.

Now, we summarize some definition and results of almost disjoint families that will be helpful in our constructions. 

\begin{mydef}
If $\cA$ is an almost disjoint family, we denote the ideal generated by $\cA$ as $\cI(\cA)$. We always assume that our ideals contain the finite sets, then $\cI(\cA)=\langle\cA\cup[\w]^{<\w}\rangle$.  
\end{mydef}

Here, $\langle\cB\rangle$ denotes the ideal generated by $\cB$. Clearly an element of $\cI(\cA)$ is a set that can be almost covered by a finite subfamily of $\cA$.

\begin{mydef}
For an almost disjoint family $\cA$ we define $\dobleposit:=\{X\subseteq \omega\mid\{A\in\cA\mid |X\cap A|=\omega\}\geq\omega\}$.
\end{mydef}

Note that
$\dobleposit\subseteq\posit$ and $\cA$ is mad if and only if, $\dobleposit=\posit$.

If $\cA$ is an almost disjoint family on a countable set $N$, then its \textit{Isbell-Mr\'owka space}, $\Psi(\cA)$,
is $N \cup \cA$, where $N$ is the set of isolated points and the basic neighborhoods of each $A\in\cA$ are of the form $\{A\}\cup A\setminus F$ for $F\in[N]^{<\omega}$. The space $\Psi(\cA)$ has interesting properties, it is Hausdorff, zero-dimensional, separable and locally compact. 
The Franklin space of $\cA$ is then the one point compactification of $\Psi(\cA)$, i.e., $\cF(\cA)=\Psi(\cA)\cup\{\infty\}$. The basic neighborhoods of $\infty$ are of the form:
\[\{\infty\} \cup (\cA\setminus\{A_0,\dots A_n\})\cup \Big(N\setminus \bigcup_{m\leq n}(A_m\cup F)\Big),
\]
where $n\in\omega$ , $A_0,\dots,A_n\in \cA$ and $F\in[N]^{<\omega}$. To learn more on Isbell-Mr\'owka spaces the reader may consult the survey \cite{topologyofmrowkaisbellspaces} or the articles \cite{Qsetsandnormality}, \cite{pseudocompactnessof}, \cite{smallmadfamilies}, \cite{maximalalmostdisjointfamiliesandpseudocompactnessofhyperspaces}, \cite{weaknormality} and \cite{madnessandweak}.

Since their early beginnings the Mr\'owka-Isbell spaces $\Psi(\cA)$ have turned out to be a fruitful source of examples and counterexamples in topology, since combinatorial properties of $\cA$ translate in topological properties of $\Psi(\cA)$. Plenty of applications of Mrowka-Isbell spaces to topology and functional analysis can be found in the literature. 
However, constructing mad families with special combinatorial properties is usually a difficult task without extra assumptions beyond $\mathsf{ZFC}$. An example of ad families with special combinatorial properties are the completely separable ones.

\begin{mydef}
An almost disjoint family $\cA$ is completely separable if for all $X\in\dobleposit$ there exists $A\in \cA$ such that $A\subseteq X$.
\end{mydef}

These families were introduced by Hechler in \cite{hechlercompletamenteseparable} and in \cite{shelaherdos}. Erd\H{o}s and Shelah asked if there exist completely separable families that are mad. Most of the early work on this question was done by Balcar and Simon in \cite{balcarsimon}, who proved that completely separable mad families exist under many assumptions. The next major advance on this was done by Shelah in \cite{shelah}, who introduced a technique that allowed him to prove that there are completely separable families under $\s<\a$ and also under $\s\geq\a$ + a certain “PCF hypothesis”. The last big progress in the search of these families was done by Mildenberger, Raghavan and Stepr\={a}ns in \cite{dilip} and \cite{onweaklytight} who improved and extended the Shelah's technique in order to get the existence of a completely separable mad family under $\s\leq\a$. This technique has turned out to be very powerful not only in the construction of completely separable mad families but also for other special mad families, e.g, weakly tight mad families (under $\s\leq\b$) and $+$-Ramsey mad families (see \cite{onweaklytight} and \cite{plusramsey} respectively). This technique will play a crucial role in Section \ref{section b=s}. The interested reader may consult \cite{hruvsak2014almost} and \cite{guzmancompletely}.

Finally, Simon proved the following remarkable theorem: 
\begin{theorem}\label{completamenteseparabledeZFC}
\cite{galvin2007vcech}
 There is a completely separable nowhere mad family.
\end{theorem}

So we can find a completely separable almost disjoint family in $\mathsf{ZFC}$ if we do not demand maximality on it. As we will see in the Section \ref{sectionseqcompactnot2ramsey}, this family will lead to a sequentially compact space not $2$-\RRamsey. The recursive construction can be carried on thanks to the following result:

\begin{lemma}{See \cite{galvin2007vcech} and \cite{guzmancompletely}.}\label{si te comes uno te comes muchos} If $\cA$ is a completely separable almost disjoint family and $X\in \dobleposit$, then the set $\{A\in\cA\mid A\subseteq X\}$ has cardinality $\mathfrak{c}$.
\end{lemma}

Recall that a space $X$ is Fr\'echet if whenever $x\in \overline{A}$, with $A\subseteq X$, there exists $\{x_n:n\in\omega\}\subseteq A$ that converges to $x$. Obviously, every first countable space is Fr\'echet, but there are many Fr\'echet spaces that are not first countable. An almost disjoint family $\cA$ is nowhere mad if for every $X\in\cI^+(\cA)$ there exists $B\in[X]^\omega$ that is almost disjoint with every element in $\cA$.

\begin{theorem}\cite{hruvsak2014almost}\label{nowhere mad sii frechet}
Let $\cA$ be an almost disjoint family, then $\cA$ is nowhere mad if and only if $\cF(\cA)$ is a Fr\'echet space.
\end{theorem}


We now remind some important results concerning the $n$-sequential compactness property proved in \cite{kubis2021topological}.

\begin{proposition}\label{prop10}
\cite{kubis2021topological} The Franklin space of a mad family is not $2$-\RRamsey.
\end{proposition}

Let $\cA$ be an almost disjoint family on a countable set $N$, $f:[\w]^n\rightarrow\cF(\cA)$ and suppose we want to see if we can find a convergent subsequence of $f$. By passing to an infinite subset of $\w$ if needed (this can be done by applying Ramsey's theorem), we may assume that one of the following cases holds:
\begin{enumerate}
    \item $f$ is constant $\infty$,
    \item $f''[\w]^n\subseteq\cA$,
    \item $f''[\w]^n\subseteq N$.
\end{enumerate}
The following lemma tells us that we only need to worry in the last case:
\begin{lemma}\label{lemma10}
\cite{kubis2021topological} If $\cA$ is an almost disjoint family on a countable set $N$ of isolated points and  $n\geq 1$, $\cF(\cA)$ is $n$-\Ramsey if and only if for all $f:[\omega]^n\rightarrow N$, there exists $B\in[\N]^\N$ such that $f\restriction [B]^n$ converges.
\end{lemma}

Recall that the \emph{character of} $x\in X$ is the minimum cardinality of a local base of $x$ and it is denoted by $\chi(x,X)$. The character of $X$ is $\chi(X)=\sup\{\chi(x,X):x\in X\}$.

\begin{theorem}\label{teoremadepaulykubis}
\cite{kubis2021topological} 
Every sequentially compact space with character less than $\b$ is $n$-\Ramsey for all $n\geq 1$.
\end{theorem}

In Theorem \ref{lashipotesisdelteodekubisypaulsonoptimas} we will see that $\b$ is optimal for the Theorem above.
A consequence of Theorem \ref{teoremadepaulykubis} and Lemma \ref{lemma10} is:

\begin{corollary}{\textnormal{\cite{kubis2021topological}}}\label{closure}
Let $\cA$ be an almost disjoint family on a countable set $N$, such that for every $f:[\omega]^n\to N$ there exists $X\in[\omega]^\omega$ so that 
$f''[X]^n$ has closure of size less than $\mathfrak{b}$ in $\Psi(\cA)$. Then $\cF(\cA)$ is $n$-\RRamsey.
\end{corollary}

We also get:

\begin{corollary}\label{lema de paul}
\cite{kubis2021topological}
Let $\cA$ be an almost disjoint family on a countable set $N$ of isolated points and $n\geq 1$. If for all $f:[\omega]^n\rightarrow N$, there exists $B\in[\N]^\N$ such that
\[
\{A\in\cA\mid |A\cap f\restriction[B]^n|=\w\}
\]
has size less than $\mathfrak{b}$, then $\cF(\cA)$ is $n$-\RRamsey.
\end{corollary}

One of the key features of the counterexamples given in \cite{kubis2021topological} is the structure of the $\textsf{FIN}^n$ ideal on $\omega^n$. This ideals were introduced by Katetov in \cite{ondescriptiveclassification}. Recall that $\textsf{FIN}=[\w]^{<\w}$ denotes the ideal of finite sets of $\w$.
The ideal $\textsf{FIN}^n$ is defined recursively as the Fubini product of $\textsf{FIN}\otimes \textsf{FIN}^{n-1}$. 
That is, $X\in\textsf{FIN}^n$ if 
$$\{x\in\w^{n-1}:|\{n\in\w:x\conc n\in X\}|=\omega\}\in\textsf{FIN}^{n-1}.$$

To learn more about this ideals, the reader may consult \cite{whenthekatetovorder}, \cite{somestructuralaspects}, \cite{onaconjectureofdebs} and \cite{ideallimitsofseq}.

The following combinatorial lemma will be very useful for us in the future.
\begin{lemma}\label{existe un infinito con imagen en fin}
\cite{kubis2021topological} For every $f:[\N]^{n}\rightarrow\N^{n+1}$, there exists $X\in[\N]^\N$ such that \emph{$f''[X]^n\in\textsf{FIN}^{n+1}$}.
\end{lemma}

To deal with these kind of functions,
 we will also need an useful tool known as the canonical Ramsey theorem. A very deep and powerful fact due to Erd\"os and Rado is that there is a finite, very easy to describe ''base'' for the equivalence relations on $[\w]^n$. 

For the rest of this paper, given $a\in[\omega]^n$, we write it as $\{a_i:i<n\}$ where $a_0<\cdots<a_{n-1}$. If $I\subseteq n$, let $a\rest I:=\{a_i:i\in I\}$. For every $I\subseteq n$, we can define the equivalence relation $E_I$ on $[\omega]^n$ given as follows:
\[aE_Ib\Longleftrightarrow a\rest I=b\rest I.\]

Note that:
\begin{itemize}
    \item $aE_\emptyset b$ for every $a,b\in[\omega]^n$ and 
    \item $aE_nb$ if and only if $a=b$.
\end{itemize}

Of course not every equivalence relation is of this form. However, the surprising fact is that every equivalence relation is ``locally'' one of the ones we just described. Formally, we have the following:

\begin{theorem}\label{erdosradoteo}{\textnormal{\cite{erdos1950combinatorial}}}
For every equivalence relation $E$ on $[\omega]^k$ there are $M\in[\omega]^\omega$ and $I\subseteq k$ such that $aEb$ if and only if $a\rest I=b\rest I$.
\end{theorem}

The reader can learn more in this topic from Todorcevic's book \cite{Todorcevic+2010}.

Applying Theorem \ref{erdosradoteo} several times, we get the following: 
\begin{corollary}\label{corodeerdosrado}
    Let $E_1,\dots, E_l$ equivalence relations on $[\w]^n$. There exist $M\in[\w]^\w$ and $I_1,\dots, I_l\subseteq n$ such that for all $a,b\in [M]^n$ and $i\leq l$, we have that \[aE_ib\Longleftrightarrow a\restriction I_i=b\restriction I_i.\]
\end{corollary}

\section{Examples of $n$-\Ramsey spaces}

In this section we will give some examples of $n$-\Ramsey spaces. By Theorem \ref{teoremadepaulykubis}, we know that there are sequentially compact spaces that fail to even satisfy the $2$-sequential compactness property and that, in general, the properties of $n$-sequential compactness and $(n+1)$-sequential compactness may be consistently different. To see how subtle is to stratify the classes of $n$-\Ramsey spaces, we will analyze some classical examples of sequentially compact non first countable spaces. These are necessary assumptions due to the following results:

\begin{proposition}\label{prop15}\cite{kubis2021topological}
    \begin{enumerate}
        \item If $n\geq 1$ and $X$ is $(n+1)$-\RRamsey, then $X$ is $n$-\RRamsey. In particular if $X$ is $n$-\Ramsey for some $n\geq 1$, then $X$ is sequentially compact.

        \item If $X$ is sequentially compact and first countable, then $X$ is $n$-\Ramsey for all $n\geq 1$.
    \end{enumerate}
\end{proposition}

The following five examples are all the sequentially compact, non first countable spaces that appear in the book \cite{contraejemplos} (which can also be conveniently accessed in the website \cite{The_pi-Base_Community_pi-Base_data}):
\begin{enumerate}
    \item The Fort space.
    \item The closed long ray.
    \item The altered long ray.
    \item The Alexandroff square.
    \item $\w_1+1$.
\end{enumerate}

We will describe them in detail below and see that indeed all of them are $n$-\Ramsey for every $n\geq 1$.
Note that as an easy consequence of Proposition \ref{prop15} item (2) we have that:
\begin{corollary}\label{corolariodelteodepaulykubis}
\cite{kubis2021topological}   Let $X$ be sequentially compact space in which the closure of every countable set is first countable. Then $X$ is $n$-\Ramsey for all $n\geq 1$.
\end{corollary}

From the previous corollary, we have that the one-point compactification of every discrete space is $n$-\Ramsey for all $n\geq 1$, in particular, the \emph{Fort space}, that is, the one-point compactification of the discrete space of size $\c$, is $n$-\Ramsey for all $n\geq 1$.

\begin{example}\label{longclosedrayspace}
(\textbf{Closed long ray}) Consider $\mathbb{L}=\w_1\times [0,1)$ with the lexicographic order and $\Tilde{\mathbb{L}}=\mathbb{L}\cup\{\infty\}$ where $\infty$ is the maximum point. Then \emph{the closed long ray} $\Tilde{\mathbb{L}}$ with the order topology is $n$-\Ramsey for all $n\geq1$.
\end{example}

\begin{proof}
Let $f:[\w]^n\rightarrow \Tilde{\mathbb{L}}$, applying Ramsey's theorem we know that there exist $B\in[\w]^\w$ such that either, $f''[B]^n\subseteq\mathbb{L}$ or $f''[B]^n=\{\infty\}$. In later case we are done and in the former note that $f''[B]^n\subseteq\mathbb{L}$ and $\mathbb{L}$ is first countable and sequentially compact so we can apply Corollary \ref{corolariodelteodepaulykubis}. 
\end{proof}

Note that $\Tilde{\mathbb{L}}$ in Example \ref{longclosedrayspace} is sequentially compact but it is not first countable, since $\infty\in \Tilde{\mathbb{L}}$ has no countable base.
\begin{example}
    \textbf{(Altered long ray)} Consider $\widehat{\mathbb{L}}=\w_1\times[0,1) \cup \{\infty\}$ where $\w_1\times[0,1)$ has the lexicographic order and the basic neighborhoods of $\infty$ are of the form:
    \[
    N(\infty,\beta)=\{\infty\}\cup \bigcup_{\alpha\geq\beta}(\{\alpha\}\times (0,1))=\{\infty\}\cup ((\w_1\setminus\beta)\times(0,1)).
    \]
    Then $\widehat{\mathbb{L}}$  is $n$-\Ramsey for all $n\geq1$.
\end{example}
\begin{proof}
    Let $f:[\w]^n\rightarrow \widehat{\mathbb{L}}$, applying Ramsey's theorem we know that there exist $B\in[\w]^\w$ such that either, $f''[B]^n\subseteq\mathbb{L}$ or $f''[B]^n=\{\infty\}$. But note that the topology of $\mathbb{L}$ as a subspace of $\widehat{\mathbb{L}}$ coincides with the topology of $\mathbb{L}$ as a subspace of $\Tilde{\mathbb{L}}$.
\end{proof}

Let us denote by $\Delta$ the set of points of the form $(x,x)$ for $x\in[0,1]$. If $x\in [0,1]$ then denote by $V_x$ the vertical line given by $x$, i.e., $V_x=\{(x,y)\mid y\in[0,1]\}$. 

The following example illustrates the usual way in which we will be applying Theorem \ref{erdosradoteo} to derive $n$-sequential compactness of certain spaces.

\begin{example}
    (\textbf{Alexandroff square}) Let $\mathbb{A}=[0,1]\times[0,1]$ with the following topology: if $p=(x,x)\in\Delta$ then a basic neighborhood of $p$ is of the form:
    \[
    N(p,x_0,\dots,x_{n-1},\epsilon)=\left([0,1]\times(x-\epsilon,x+\epsilon)\right)\setminus(\bigcup_{i\in n}V_{x_i}),
    \]
    where $n\in\w$, $\epsilon>0$ and $x\not\in\{x_i\mid i\in n\}$. If $p=(x,y)\not\in\Delta$ then a basic neighborhood of $p$ is of the form $N(p,\epsilon)=\{(x,t)\mid |y-t|<\epsilon\}$. Then $\mathbb{A}$ is $n$-\Ramsey for all $n\geq 1$. 
\end{example}
\begin{proof}
    Let $f:[\w]^n\rightarrow \mathbb{A}$. Now for every $i\in 2$ let $f_i:[\w]^n\rightarrow[0,1]$ such that for all $s\in[\w]^n$ we have that $f(s)=(f_0(s),f_1(s))$.
    Note that: 
    \begin{enumerate}[label=(\roman*)]
        \item\label{2defacts} If there is $B\in[\w]^\w$ such that $f_0\restriction[B]^n$ is constant, then $f\restriction [B]^n$ has a convergent subsequence\footnote{In fact, more can be said here: if there is $p\in \mathbb{A}\setminus\Delta$ and $B\in[\w]^\w$ such that $f\restriction[B]^n$ converges to $p$, then $f_0$ is constant on $[B]^n$.}.
    \end{enumerate}

    Applying Theorem \ref{erdosradoteo} to the function $f_0$, we can get $M\in[\w]^\w$ and $I\subseteq n$ such that for all $s,t\in[M]^n$, $f_0(s)=f_0(t)$ if and only if $s\restriction I=t\restriction I$. 
    
    If $I=\emptyset$ then we are done by \ref{2defacts}. So assume that $I\neq\emptyset$ and let us consider $f_1\restriction[M]^n$. As $[0,1]$ is $n$-\Ramsey (with the usual topology), there exists $N\in[M]^\w$ and $x\in[0,1]$ such that $f_1\restriction[N]^n$ converges to $x$. We claim that $f\restriction [N]^n$ converges to $p=(x,x)$. So let $N(p,x_0,\dots,x_{n-1},\epsilon)$ be a basic neighborhood of $p$. As $f_1\restriction [N]^n$ converges to $x$, then there exists $F\in [N]^{<\w}$ such that $f_1''[N\setminus F]^n\subseteq(x-\epsilon,x+\epsilon)$. 
    
    Now, as $I\neq\emptyset$, we claim that there is $G\in[N]^{<\w}$  such that $f_0''[N\setminus G]^n\cap \{x_i\mid i\in n\}=\emptyset$. Indeed, for every $i\in n$ there is $s\in [\w]^{<\w}$ not empty such that $f_0(t)=x_i$ if and only if $t\restriction I=s_i$. Now, $G=\bigcup_{i\in n}s_i$ is as desired.

    This shows that $f''[N\setminus(F\cup G)]^n\subseteq N(p,x_0,\dots,x_{n-1},\epsilon)$ and then $f\restriction [N]^n$ converges to $p$.
\end{proof}

Finally, we will prove that in the realm of limit ordinals, being sequentially compact, having uncountable cofinality and being $n$-\Ramsey for all $n\geq 1$ are equivalent. We will first need the following:

\begin{proposition}\label{cerradura numerable}
    If $\gamma$ is an ordinal and $X\subseteq\gamma$ is countable, then $\overline{X}$ is countable.
\end{proposition}
\begin{proof}
    Suppose that $\overline{X}$ is not countable and let $Y$ of size $\w_1$ such that $X\subseteq Y\subseteq\overline{X}$. As $Y$ is a set of ordinals, by going to a subset if necessary, we can enumerate $Y$ in order type $\w_1$, i.e., $Y=\{x_\alpha\mid\alpha\in\w_1\}$ where $x_\alpha<x_\beta$ if $\alpha<\beta$. Note that as $X$ is countable, there exists $\alpha\in\w_1$ such that $\{x_\beta\mid \alpha<\beta<\w_1\}\cap X=\emptyset$, so in particular $(x_\alpha,x_\beta]\cap X=\emptyset$ for any $\beta>\alpha$, which shows that $x_\beta\not\in \overline{X}$, a contradiction.
\end{proof}

Note that by the last proposition, if $\gamma$ is a limit ordinal with uncountable cofinality and $X\subseteq \gamma$ is countable, then the closure of $X$ in $\gamma$ is sequentially compact and first countable, so by Corollary \ref{corolariodelteodepaulykubis} we have that:



\begin{corollary}\label{los ordinales son n ramsey si y solo si son seq comp}
Let $\gamma$ be a limit ordinal. Then $cof(\gamma)>\w$ if and only if $\gamma$ with the order topology is $n$-\Ramsey for all $n\geq 1$. 
\end{corollary}


Finally, we get a full characterization of those ordinals that are $n$-\Ramsey for all $n\in\omega$.

\begin{corollary}
    An ordinal $\gamma$ is $n$-\Ramsey for all $n\geq 1$ if and only if $cof(\gamma)\neq\omega$.
\end{corollary}

\begin{proof}
Let $f:[\omega]^n\to \gamma$ and assume that $\gamma$ has uncountable cofinality. Let $X=f''[\omega]^n$ and define $\alpha=\sup(X)$. Since $\gamma$ does not have countable cofinality, $\alpha\in\gamma$. It follows then that $\alpha=\max(\overline{X})$. Take $\overline{X}$ with the inherit order, since the maximum of $\overline{X}$ exists, $\beta=ot(\overline{X})$ is a countable successor ordinal. Therefore it is first countable and sequentially compact, which implies that $\beta$ is $n$-\Ramsey by Proposition \ref{prop15} item (2). In consequence, $f$ has a convergent subsequence in $\overline{X}\subseteq\gamma$.
\end{proof}


With this we now know that all sequentially compact examples in \cite{contraejemplos} (and only those) are $n$-\Ramsey for all $n\geq 1$.

\section{Sequentially compact Fr\'echet spaces that are not $2$-\Ramsey}\label{sectionseqcompactnot2ramsey}

It is well know that in the case of Franklin spaces sequentially compact does not imply being Fr\'echet, however, Paul Szeptycki pointed out that a $2$-\Ramsey Franklin space must be Fr\'echet.


\begin{proposition}\label{2ramseyimplicafrechet}    
If $\cA$ is an almost disjoint family such that $\cF(\cA)$ is $2$-\RRamsey, then $\cF(\cA)$ is Fr\'echet.
\end{proposition}
\begin{proof}
    Suppose that $\cF(\cA)$ is not Fr\'echet, this means, by Theorem \ref{nowhere mad sii frechet}, that $\cA$ is somewhere mad, i.e., there exists $X\in\cI^+(\cA)$ such that $\cA\restriction X:=\{A\cap X\mid A\in\cA(|A\cap X|=\w)\}$ is an infinite mad family on $X$ and then $\cF(\cA\restriction X)$ is not $2$-\Ramsey by Proposition \ref{prop10}.

    Now let $Y:=X\cup\{A\mid A\in\cA\restriction X\}\cup\{\infty\}$ with the topology inherited by $\cF(\cA)$. It is easy to see that $Y$ is a closed subset of $\cF(\cA)$ 
    and it is homeomorphic to $\cF(\cA\restriction X)$. As $\cF(\cA)$ is $2$-\Ramsey and $Y\subseteq\cF(\cA)$ is closed, $Y$ is $2$-\Ramsey and consequently $\cF(\cA\restriction X)$ is $2$-sequentially compact, which is a contradiction.
    \end{proof}

In \cite{kubis2021topological}, the authors asked whether there is a \textsf{ZFC} example of a Fr\'echet sequentially compact space that is not $2$-\RRamsey. They showed that a completely separable mad family suffices for the construction of such example. As discussed earlier, the existence of a completely separable mad family holds under very weak assumptions like $\c<\aleph_\w$, however, as we stated in Theorem \ref{completamenteseparabledeZFC}, completely separable almost disjoint families, possibly not maximal, do exist in \textsf{ZFC}. The existence of one of these families in $\textsf{ZFC}$ is good enough for the construction sketched in \cite{kubis2021topological}, we give a full proof below.

If $N$ is countable set, $\cA\subseteq\cP(N)$ and $X\in[N]^\w$, then  $X\perp\cA$ means that $X\cap A$ is finite for all $A\in\cA$. The family $\{X\in[N]^\w\mid X\perp\cA\}$ is denoted by $\cA^\perp$. .

\begin{lemma}\label{noramsey}
    Let $\cA$ an almost disjoint family on a countable set $N$ of isolated points. If there exists $\cA_0=\{A_m\mid m\in\w\}\subseteq[N]^\omega$ and $f:[\w]^n\rightarrow N$ such that:
    \begin{enumerate}[label=(\roman*)]
        \item\label{lemanorramsey0} Either, $\cA_0\subseteq\cA$ or $\cA_0\subseteq\cA^\perp$.
        \item\label{lemanorramsey1} $\cA_0$ is a pairwise disjoint family. 
        \item\label{lemanorramsey2} $f$ is one-to-one.
        \item\label{lemanorramsey3} If $\min(s)=m$ then $f(s)\in A_m$.
        \item\label{lemanorramsey4} For all $B\in [\w]^\w$ there exists $A\in\cA\setminus\cA_0$ such that $|f''[B]^n\cap A|=\w$.
    \end{enumerate}
    Then $\cF(\cA)$ is not $n$-\RRamsey.
\end{lemma}
\begin{proof}
We will show that $f\restriction [B]^n$ does not converge to any point in $\cF(\cA)$ for all $B\in[\w]^{\omega}$.
Fix $B\in[\w]^{\omega}$, then as $f\restriction [B]^n$ is one-to-one, it does not converge to any isolated point. Let $F\in[B]^{<\w}$ and $m\in\omega$. There are $m'>m$ and $a\in[B\setminus F]^n$ such that $\min(a)=m'$. Thus $f(a)\in A_{m'}$. This shows that $f''[B\setminus F]^n\nsubseteq \{A_m\}\cup A_m$ for any $F$ and thus $f\rest[B]^n$ does not converge to $A_m$ in the case where $\cA_0\subseteq\cA$.

Now let us see that $f\restriction [B]^n$ does not converge to any point in $\cA\setminus\cA_0$ either. Let $F\in[B]^{<\w}$ and fix $A\in\cA\setminus\cA_0$. Take $m=\min(B\setminus F)$. Thus by condition \ref{lemanorramsey3}, $f''[B\setminus F]^n\cap A_m$ is infinite and since $A\cap A_m$ is finite, it follows that $f''[B\setminus F]^n\nsubseteq \{A\}\cup A$.

The only thing that remains to be seen is that $f\restriction [B]^n$ does not converge to $\infty$. Let $A\in\cA\setminus\cA_0$ given by item \ref{lemanorramsey4}. We will show that $f''[B\setminus F]^n\nsubseteq U:=\cF(\cA)\setminus (\{A\}\cup A)$ for all $F\in [B]^{<\w}$. Since $U$ is a basic neighborhood of $\infty$ in $\cF(\cA)$, this will suffice. Let $m=\max F$. Since $A\cap A_i$ is finite for every $i\leq m$ and $|f''[B]^n\cap A|=\w$, there exists $a\in [B\setminus F]^n$ such that $f(a)\in A$ and we are done.

This shows that $\cF(\cA)$ is not $n$-\Ramsey as witnessed by $f$.
\end{proof}

We can now prove the following, answering a question in \cite{kubis2021topological}.

\begin{theorem}
There exists (in $\mathsf{ZFC}$) a space which is Fr\'echet, sequentially compact and it is not $2$-sequentially compact.
\end{theorem}
\begin{proof}
Let $\cA$ be an almost disjoint family on $\omega\times\omega$ which is completely separable, infinite and nowhere mad. We can assume that $A_n=\{n\}\times\omega\in\cA$ for all $n\in\omega$.
Now let $G:[\omega]^2 \rightarrow \cF(\cA)$ given by $G(\{m,n\})=(m,n)$, where $m<n$. We will show that $\cA$, $\cA_0=\{A_n\mid n\in\w\}$ and $f=G$ satisfy the hypothesis of Lemma \ref{noramsey} for $n=2$.

Conditions \ref{lemanorramsey1}, \ref{lemanorramsey2} and \ref{lemanorramsey3} are clear. For \ref{lemanorramsey4} note that \ref{lemanorramsey2} and \ref{lemanorramsey3} imply that $|G''[B]\cap A_m|=\w$ for every $m\in B$, in particular $G''[B]^2\in\dobleposit$. Since $\cA$ is completely separable, using Lemma \ref{si te comes uno te comes muchos}, there exists $A\in\cA\setminus\cA_0$ as required in \ref{lemanorramsey4}.
\end{proof}

\section{$n$-\Ramsey spaces from $\mathfrak{b}=\c$.}

The relevance of the cardinals $\b$ and $\d$ in our examples is evident due to the next few facts regarding $\textsf{cov}^*(\textsf{FIN}^n)$ and $\textsf{cof}(\textsf{FIN}^n)$. We will first introduce the concept of Hechler trees in order to make the analysis of these results cleaner.

Given $n\in\omega$, we will say that $T\subseteq\omega^{\leq n+1}$ is a {\it Hechler tree}, if for every $s\in T\cap\omega^{\leq n}$, the set of successors of $s$ is \emph{coinitial} on $\w$, i.e., there exists $k_s\in\omega$ such that $succ_T(s)=\{i\in\omega\mid s\conc i\in T\}=\omega\setminus k_s$. We can naturally associate a function $h_T:T\cap\omega^{\leq n}\to\omega$ to each Hechler tree $T\subseteq\omega^{\leq n+1}$ by defining $h(s)=k_s$. 

For $s,t\in\omega^{\leq n}$, we say that $s\prec t$ if either $s\sqsubset t$ or $s=r\conc l$, $t=r\conc m$ for some $r\in\omega^{\leq n}$ and $l<m$, where $s\sqsubset t$ means that $s$ is a proper initial segment of $t$.
Enumerate $\omega^{\leq n}=\{s_i\mid i\in\omega\}$ such that $s_i\prec s_j$ implies $i<j$. We fix this enumeration once and for all.

\begin{mydef}
    Let $f\in\w^\w$ and let $T\subseteq\w^{n+1}$ be a Hechler tree. We define
    \begin{itemize}
        \item $T_f\subseteq\omega^{\leq n+1}$ where $\emptyset\in T_f$ and $succ(s_i)=\omega\setminus f(i)$ for every $s_i\in T_f\cap\omega^{\leq n}$.
        \item $f_T\in\omega^\omega$ where $f_T(n)=k$ if and only if $s_n\in T$ and $succ_T(s_n)=\omega\setminus k$ and $f_T(n)=0$ otherwise.
    \end{itemize}
\end{mydef}
Let $C_s=\{t\in\omega^{n+1}\mid s\sqsubseteq t\}$ for every $s\in\w^{\leq n}$ and let us write $C_k$ instead of $C_{(k)}$.
Also, for a tree $T\subseteq\omega^{\leq n}$ we denote its set of branches by $[T]=T\cap\omega^n$.
It follows directly from the definition that for every $X\in\textsf{FIN}^n$ there exists a Hechler tree $T\subseteq\w^{\leq n}$ such that $X\cap [T]=\emptyset$.

\begin{lemma}\label{cov igual a b}
For every $n\geq 1$, if \emph{$\cH\subseteq \textsf{FIN}^{n+1}$} is such that $|\cH|<\b$ and $B\in [\omega]^\w$, then there exists a block sequence\footnote{Recall that $(a_m\mid m\in\w)\subseteq \w^n$ is a block sequence if $\max a_m<\min a_{m+1}$ for all $m\in\w$ and each $a_m$ is increasing.} $A\subseteq B^{n+1}$ such that $|A\cap X|<\w$ for every $X\in\cH$.
\end{lemma}

\begin{proof}
For every $X\in\cH$ let $T_X$ be a Hechler tree disjoint from $X$. We write $f_X$ instead of $f_{T_X}$.

Since $|\cH|<\b$, there is a function $f\in\w^\w$ such that $f>^*f_X$ for all $X\in\cH$. We will recursively define a block sequence $(a_m\mid m\in\omega)\subseteq\w^{n+1}$. Let $(a_m^i:i\leq n)$ denote the sequence $a_m$. Assume we have constructed $a_k$ for every $k<m$. For $a_m$ define $a_m^0\in B\setminus\left(f(\emptyset)\cup\max(a_{m-1})\right)$ (where $a_{-1}=0$). Now for every $i<n$ pick
$$a_m^{i+1}\in B\setminus f((a_m^0,\ldots,a_m^i))$$
bigger that $a_m^i$. These choices are possible since $B$ is infinite. This finishes the construction of $A=(a_m\mid m\in\omega)\subseteq B^{n+1}$. 

To see that it works fix $X\in\cH$. Since $f>^*f_X$ there exists $k\in\omega$ such that if $f_X(n)>f(n)$ then $s_n\in\bigcup_{i<k}C_i$. As $(a_m\mid m\in\w)$ is a block sequence, $\{a_m^0\mid m\in\w\}$ is increasing and thus $a_m^0>\max\{k,f_X(\emptyset)\}$ for all but finitely many $m\in\w$. 
Given any $a_m^0$ as above, $a_m^0>f_X(\emptyset)$ and as $(a_m^0,\ldots,a_m^i)$ extends $a_m^0$ in $\w^{n+1}$ for any $i\leq n$, it follows from the assumption over $k$ that $f((a_m^0,\ldots,a_m^i))>f_X((a_m^0,\ldots,a_m^i))$. In particular $(a_m^0,\ldots,a_m^n)\in T_X$ which implies that $X\cap A\subseteq \{a_m\mid a_m^0<k\}$, with the latter set being finite.
\end{proof}

The previous lemma gives in particular that $\b\leq\textsf{cov}^*(\textsf{FIN}^{n+1})$ for every $n\geq 1$. The other inequality can also be proved using similar ideas.

\begin{proposition}
\emph{$\textsf{cov}^*(\textsf{FIN}^{n+1})=\b$} for every $n\geq 1$.
\end{proposition}

\begin{proof}
We only need to prove that $\textsf{cov}^*(\textsf{FIN}^{n+1})\leq\b$. Suppose that $\kappa<\textsf{cov}^*(\textsf{FIN}^{n+1})$ and let $\{h_\beta\in\w^\w\mid \beta\in\kappa\}$. Recall that $\b$ is the minimum size of a family $\cB\subseteq\w^\w$ that is unbounded on each infinite subset of $\w$ (see \cite{van1984integers}), so it is enough to see $\{h_\beta\in\w^\w\mid \beta\in\kappa\}$ is bounded on some $B\in[\w]^\w$. Enumerate $\omega^n$ as $\{t_l\mid l\in\omega\}$ such that if there is $r\in\omega^{n-1}$ such that $t_l=r\conc i$ and $t_{l'}=r\conc j$ and $i<j$ then $l<l'$.

    For each $\beta\in\kappa$ consider the set $Y_\beta=\{t_l\conc m\mid l\in\w \wedge m\leq  h_\beta(l)\}$. Clearly each $Y_\beta$ is in $\textsf{FIN}^{n+1}$. Now consider $\cH=\{Y_\beta\mid \beta\in\kappa\}\cup\{C_{t_l}\mid l\in\omega\}$. It is clear that $|\cH|<\textsf{cov}^*(\textsf{FIN}^{n+1})$, so there exists $A\in [\omega^{n+1}]^\w$ such that $A\cap C_{t_l}=^*\emptyset$ for all $l\in\w$ and $A\cap Y_\beta=^*\emptyset$ for every $\beta\in\kappa$.

    Now define $f:\w\rightarrow\w$ as follows: 
    \begin{center}
            $f(l)= \left\{
    	       \begin{array}{ll}
                    \max\{i\in\omega\mid t_l\conc i\in A\}       & \mathrm{if\ } A\cap C_{t_l}\neq\emptyset \\
                    0      & \text{ otherwise.}
    	       \end{array}
	     \right.$
        \end{center}
Let $B=\{l\in\omega:A\cap C_{t_l}\neq\emptyset\}\in[\omega]^\w$.
Now, as $A\cap Y_\beta=^*\emptyset$, there exists $l\in\w$ such that if $m>l$ then $C_{t_m} \cap (A\cap Y_\beta)=\emptyset$, which means that $f(m)>h_\beta(m)$ whenever $m\in B\setminus (l+1)$, i.e., $f$ dominates $h_\beta$ on $B$.
\end{proof}


Now we are ready to prove that under $\mathfrak{b}=\mathfrak{c}$ we can distinguish the classes of $n$-\Ramsey and $(n+1)$-\Ramsey spaces for every $n\geq 1$.

\begin{theorem}
\textnormal{($\mathfrak{b}=\mathfrak{c}$)} For all $n\geq 2$ there exists a $n$-\Ramsey space which is not $(n+1)$-\RRamsey.
\end{theorem}
\begin{proof}
We will construct an almost disjoint family $\cA$ 
on $\N^{n+1}$ such that its Franklin space is $n$-\Ramsey but not $(n+1)$-sequentially compact. Let $G:[\omega]^{n+1} \rightarrow \cF(\cA)$ given by $G(\{m_0,\dots,m_{n}\})=(m_0,\dots,m_{n})$, where $m_0<\dots<m_{n}$.


Let us enumerate $(\N^{n+1})^{[\N]^n}=\{f_\alpha\mid\N\leq\alpha<\cont\}$ and $[\N]^\N=\{B_\alpha\mid\N\leq\alpha<\cont\}$.
We will construct $\{A_\alpha\mid\alpha<\cont\}\subseteq\cP(\omega^{n+1})$ and $\{X_\alpha\mid \N\leq\alpha<\cont\}\subseteq\cP(\omega)$ such that the following hold:
\begin{enumerate}
    \item\label{1 de b=c} For all $n\in\N$, $A_n=C_n$ and $A_\alpha$ is a block sequence for every $\alpha\geq\w$ (hence $A_\alpha\in\textsf{FIN}^{n+1}$).
    \item\label{2 de b=c}$\forall\beta,\alpha\in\cont\ (\beta\neq\alpha\rightarrow A_\beta\cap A_\alpha=^*\emptyset )$.
    \item\label{3 de b=c} $f_\alpha''[X_\alpha]^n\in \textsf{FIN}^{n+1}$.
    \item\label{4 de b=c} $\forall\beta,\alpha\in\cont\ (\N\leq\beta<\alpha\rightarrow f_\beta''[X_\beta]^n\cap A_\alpha=^*\emptyset )$.
    \item\label{5 de b=c} $\forall\alpha\in[\w,\cont)\ (|A_\alpha\cap G''[B_\alpha]^{n+1}|=\N)$.
\end{enumerate}

Assume that we have already constructed $\{A_\beta\mid \beta<\alpha\}$ and $\{X_\beta\mid \beta<\alpha\}$ for some $\alpha\geq\omega$.
We can apply Lemma \ref{cov igual a b} to the family $\cH$ defined as follows:
$$\cH=\{f_\beta''[X_\beta]^n\mid \omega\leq\beta<\alpha\}\cup\{A_\beta\mid \beta<\alpha\}$$ 
and the infinite set $B_\alpha$. In this way we get a block sequence $A_\alpha\subseteq B_\alpha^{n+1}$ such that:
\begin{enumerate}[label=(\alph*)]
    \item $\forall\beta\in\alpha\ (|A_\beta\cap A_\alpha|<\w)$.
    \item $\forall\beta\in\alpha\ (\omega\leq\beta<\alpha\rightarrow |f_\beta''[X_\beta]^n\cap A_\alpha|<\w)$.
\end{enumerate}

This ensures that conditions (\ref{2 de b=c}) and (\ref{4 de b=c}) hold, also,  as $A_\alpha$ is a block sequence of elements of $B_\alpha$, condition (\ref{5 de b=c}) holds. Now, find $X_\alpha\in[\N]^\N$ such that $f_\alpha''[X_\alpha]^n\in\textsf{FIN}^{n+1}$ (and this one exists by Lemma \ref{existe un infinito con imagen en fin}). This finishes the construction. 

Consider $\cA=\{A_\alpha\mid\alpha\in\cont\}$. 
Condition (\ref{4 de b=c}) implies that $f_\beta''[X_\beta]^n$ has closure of size less than $\cont$ in $\cF(\cA)$ and,  as $\mathfrak{b}=\cont$, it has closure of size less than $\mathfrak{b}$, so by Corollary \ref{lema de paul}, it has a convergent subsequence. On the other hand, condition (\ref{5 de b=c}) implies that condition \ref{lemanorramsey4} of Lemma \ref{noramsey} holds and the other three conditions when $f=G$ are clear, so $\cF(\cA)$ is not $(n+1)$-\RRamsey.
\end{proof}

\section{Small $n$-\Ramsey spaces.}

Parametrized diamonds were defined by Dzamonja, Hru\v{s}\'ak and Moore in \cite{moore2004parametrized}. These principles are weak versions of the well known Jensen's diamond principle $\diamondsuit$ but they have the advantage that are consistent with $\neg\textsf{CH}$. Moreover, many constructions that follow from $\diamondsuit$ can be carried with an appropriated parametrized diamond. These weak diamond principles can be defined from abstract cardinal invariants. Following \cite{vojtas1991generalized}, we say that triple $(A,B,E)$ is an {\it invariant} if $A$ and $B$ are sets of cardinality at most continuum and $E$ is a relation between $A$ and $B$, such that every element in $A$ is $E$-related with some element in $B$ and for every $b\in B$ there is an $a\in A$ such that $(a,b)\notin E$.

The evaluation $\langle A,B,E\rangle$ of an invariant $(A,B,E)$ is defined as follows:
$$\langle A,B,E\rangle=\min\{|X|:X\subseteq B\textnormal{ and }(\forall a\in A\ \exists b\in B\ (aEb))\}.$$

An invariant $(A,B,E)$ is called a \emph{Borel invariant} \cite{blassreductions} if $A$, $B$ and $E$ are Borel subsets on some Polish space. Given a Borel invariant $(A,B,E)$, its parametrized diamond principle is defined as follows:

$$\diamondsuit(A,B,E)\equiv\ \forall F:2^{<\omega_1}\to A\ \textnormal{Borel}\ \exists g:\omega_1\to B\ \forall f\in2^{\omega_1}$$
$$\ \ \ \ \ \ \ \ \ \ \ \ \ \{\alpha\in\omega_1:F(f\rest\alpha) E g(\alpha)\}\textnormal{ is stationary}.$$

In this context, $F:2^{<\w_1}\rightarrow A$ is \emph{Borel}, if for every $\delta<\omega_1$ the restriction of $F$ to $2^\delta$ is a Borel map. 
There exists also a notion of parametrized diamond where the invariant and the function $F$ above do not have the constraint of being Borel. This notion is denoted by $\Phi(A,B,E)$, nevertheless, most of the applications of these diamond principles have been shown to follow from its definable counterpart, and unlike them, the definable ones hold in many models in which the associated cardinal invariant is small.

For well known cardinal invariants, we write its evaluation instead of its triple. For example, $\diamondsuit(\b)\equiv\diamondsuit(\w^\w,\w^\w,\geq^*)$.
Some examples of the use of parametrized diamonds in topology are the following:

\begin{enumerate}
    \item \cite{moore2004parametrized} $\diamondsuit(\b)\implies\a=\omega_1$,
    \item \cite{moore2004parametrized} $\diamondsuit(\textsf{non}(\cM))$ implies that there is a Souslin line,
    \item \cite{moore2004parametrized} $\diamondsuit(\s^\omega)$ implies that there is a perfectly normal countably compact noncompact space,
    \item \cite{moore2004parametrized} $\diamondsuit(\RR,\neq)$ implies that there are no $Q$-sets,
    \item \cite{moore2004parametrized} $\diamondsuit(\r)$ implies that there is a $P$-point of character $\omega_1$,

    \item \cite{countableirresolvable} $\diamondsuit(\r_\textsf{scat})$ implies that there is a countable $T_3$ irresolvable space of weight $\omega_1$,

    \item \cite{scatteredspaces} $\diamondsuit(\s)$ implies that there is a family of sequentially compact spaces $\{X_\alpha:\alpha<\kappa\}$ such that $\prod_{\alpha<\kappa}X_\alpha$ is not countably compact,

    \item \cite{scatteredspaces} $\diamondsuit(\b^*)$ implies the existence of a separable Jakovlev space of size $\omega_1$,

    \item \cite{scatteredspaces} $\diamondsuit(\b^*)$ implies that for every $\eta\in\omega_1$ there exists a separable compact sequential scattered space of sequential order $\eta+1$,

    \item \cite{diamantesyultrafiltrosunion} $\diamondsuit(\mathfrak{hom}_H)$ implies the existence of stable ordered-union ultrafilters of character $\omega_1$.

    \item \cite{frechetlikeproperties} $\diamondsuit(\b)$ implies that there is a Fr\'echet $\alpha_3$ space defined from an almost disjoint family that is not bisequential,

    \item \cite{morgancovering} $\diamondsuit(\omega,<)$ implies that every separable, locally compact space with property (a) has countable extent,


    \item \cite{gruffultrafilters} $\diamondsuit(\r_P)$ implies the existence of gruff ultrafilters.
    
    \item \cite{Malykhinsproblem} $\diamondsuit(2,=)$ implies that there exists a Fr\'echet and countable topological group that is not metrizable.
\end{enumerate}



We will consider partial functions in $\omega^\omega$ and they will be denoted by $f;\omega\to\omega$. Let $\mathbb{P}=\{f;\omega\to\omega:|dom(f)|=\omega\}$. Hence, we can look at the cardinal invariant defined by $(\mathbb{P},\omega^\omega,\nsucceq)$, (where we say that $f\nsucceq g$ if $\{n\in\omega:n\in dom(f)\land f(n)<g(n)\}$ is infinite). It is known that its evaluation $\langle \mathbb{P},\omega^\omega,\nsucceq\rangle$ turns out to be $\mathfrak{b}$ (see \cite{van1984integers}). However, in order to prove that its parametrized $\diamondsuit$-principles are also equivalents, we need to show that both invariants are \emph{Borel Tukey equivalents}.

\begin{mydef}\cite{blassreductions}
    Given two Borel invariants $(A,B,E)$ and $(A',B',E')$, we say that $(A,B,E)\leq^B_T(A',B',E')$ if there are Borel maps $\phi:A\to A'$ and $\psi:B\to B'$ such that $(\phi(a),b)\in E'$ implies $(a,\psi(b))\in E$.
\end{mydef}

We say that two invariants $(A,B,E)$ and $(A',B',E')$ are Borel Tukey equivalents (denoted by $(A,B,E)\equiv_T(A',B',E')$) if $(A,B,E)\geq^B_T(A',B',E')$ and $(A,B,E)\leq_T^B(A',B',E')$.
The general notion of Tukey reductions not restricted to Borel sets (often denoted by $\leq_T$) was defined in \cite{vojtas1991generalized}. Since we are only dealing with Borel invariants, we will use $\leq_T$ 
 and $\equiv_T$ for the Borel versions. 

The importance of Tukey relations between invariants in our context relies in the following result:

\begin{proposition}\cite{moore2004parametrized}\label{Tukeydiamond}
If $(A,B,E)\leq^B_T(A',B',E')$, then $\diamondsuit(A',B',E')$ implies $\diamondsuit(A,B,E)$.
\end{proposition}

We will prove that $\diamondsuit(\b)$ is equivalent to $\diamondsuit(\mathbb{P},\w^\w,\nsucceq)$ by showing that $(\w^\w,\w^\w,\ngeq^*)\equiv_T(\mathbb{P},\w^\w,\nsucceq)$. 

\begin{proposition}
    $\diamondsuit(\b)$ is equivalent to $\diamondsuit(\mathbb{P},\w^\w,\nsucceq)$.
\end{proposition}

\begin{proof}
Let us begin by noting that $(\mathbb{P},\w^\w,\nsucceq)$ is a Borel invariant. We can identify $\mathbb{P}$ as a subset of $\w^\w$ if we add $\{-1\}$ to the range of the functions and assume that $f(n)=-1$ whenever $f$ is undefined on $n$. Formally, let $W=\w\cup\{-1\}$. Clearly $W^\w\simeq\w^\w$ and 
$$\mathbb{P}=\Big\{f\in W^\w:\forall m\in\w\ \exists n>m\big(f(n)\neq -1\big)\Big\}$$
which is a Borel subset as $f(n)\neq -1$ is an open condition.
Similarly we can check that
$$f\succeq g\Leftrightarrow\forall m\in\w\ \exists n>m\ \bigg(\Big(f(n)\neq 1\Big)\land\Big(f(n)<g(n)\Big)\bigg)$$
and then $\succeq$ is a Borel relation.

To see that $(\w^\w,\w^\w,\ngeq^*)\leq_T(\mathbb{P},\w^\w,\nsucceq)$, just take $\phi:(\w^\w)\to\mathbb{P}$ as the inclusion, that is, $\phi(f)$ is the function $f$ with its codomain replaced by $W$, and $\psi:\w^\w\to\w^\w$ is the identity. It is easy to see that if $\phi(f)\nsucceq g$ then $f\ngeq^* \psi(g)$.

For the other inequality define $\phi:\mathbb{P}\to\w^\w$ such that $\phi(f)(n)=f(\hat{n})$ where $\hat{n}=\min\{m>n:f(n)\neq -1\}$. Let us also define $\psi:\w^\w\to\w^\w$ by $\psi(g)(n)=\max\{g(i):i\leq n\}$. Clearly $\psi$ is Borel. Before showing that $\phi$ es Borel, let us prove that $\phi$ and $\psi$ work.

Assume $\phi(f)\ngeq^*g$. Let $A=\{n\in\omega:g(n)>\phi(f)(n)\}$ which is an infinite set. Given $n\in A$ and $\hat{n}$ defined as above:
\[\psi(g)(\hat{n})\geq g(n)>\phi(f)(n)=f(\hat{n}).\]
Thus $\hat{A}=\{\hat{n}:n\in A\}$ is infinite and shows that $f\nsucceq\psi(g)(n)$.

It remains to show that $\phi$ is Borel. For this notice that $\{S_{n,m}:n,m\in\omega\}$, where $S_{n,m}=\{f\in\w^\w:f(n)=m\}$, is a subbase for the topology on $\w^\w$. It then suffices to see that $\phi^{-1}[S_{n,m}]$ is Borel for every $n,m\in\w$.
But this is clear since
$$f\in\phi^{-1}[S_{n,m}]\Leftrightarrow\exists k\geq m\ \big(f(k)=m\land\big(\forall j\in\omega(m\leq j<k\rightarrow f(k)=-1)\big)\big).$$

By Proposition \ref{Tukeydiamond} the corresponding diamond principles are equivalent.
\end{proof}

We are now ready to state the parametrized diamond we will use, in the right representation for our purposes:

$$\diamondsuit(\mathfrak{b})\equiv\ \forall F:2^{<\omega_1}\to\mathbb{P}\ \textnormal{Borel}\ \exists g:\omega_1\to\omega^\omega\ \forall f\in2^{\omega_1}$$
$$\ \ \ \ \ \{\alpha\in\omega_1:F(f\rest\alpha)\ngeq^*g(\alpha)\}\textnormal{ is stationary}.$$

We call the function $g$ in the previous statement, a {\it $\diamondsuit(\mathfrak{b})$-sequence}, and we will say that $g$ {\it guesses} $f$ at $\alpha$ via $F$ if $F(f\rest\alpha)\ngeq^*g(\alpha)$. We will often omit mentioning the function $F$ when no confusion arises. 

The existence of $n$-\Ramsey spaces that fail to be $(n+1)$-\Ramsey under $\mathfrak{b}=\mathfrak{c}$ can be adapted to obtain the same counterexamples with size $\omega_1$ under $\diamondsuit(\mathfrak{b})$. 
However, even though we can achieve the recursion that guarantees that the space is not $(n+1)$-\Ramsey using the guessing property of the $\diamondsuit(\mathfrak{b})$-sequence, we also have to ensure that for every sequence $f:[\omega]^n\to\omega^{n+1}$, there is a subset $X\in[\omega]^\omega$ such that the closure of $f''[X]^n$ has size less that $\mathfrak{b}$ (hence it has countable closure since $\diamondsuit(\mathfrak{b})$ implies $\mathfrak{b}=\omega_1$). We can not achieve this by using the guessing property since we will need every of these functions being guessed cofinally.
To avoid this trammel, we will make use of the fact that $\mathfrak{d}$ is the cofinality of the ideal $\textsf{FIN}^n$. Moreover, under the assumption $\mathfrak{b}=\mathfrak{d}$ we can even get a cofinal family in $\textsf{FIN}^{n}$ with stronger properties.

We will use the notation preceding Lemma \ref{cov igual a b} for the following result. In particular, we will use again Hechler trees and the same enumeration $\{s_n:n\in\omega\}$ for $\omega^{n+1}$.

\begin{lemma}\label{hechlertrees}
If $\mathfrak{b}=\mathfrak{d}$, there is a family $\{T(\alpha):\alpha<\mathfrak{d}\}$ of Hechler trees in $\omega^{\leq n+1}$ such that
\begin{enumerate}
    \item\label{hechler1} $\forall\alpha<\beta\ \exists N\in\omega\ (T(\beta)\subseteq T(\alpha)\cup\bigcup_{i<N}C_i)$ and
    \item\label{hechler2} \emph{$\forall X\in\textsf{FIN}^{n+1}\ \exists\alpha<\mathfrak{d}\ \forall\beta>\alpha\ \exists N\in\omega\ (X\cap T(\beta)\subseteq\bigcup_{i<N}C_i)$}.
\end{enumerate}
In particular there is a base for \emph{$\textsf{FIN}^{n+1}$} of size $\mathfrak{d}$.
\end{lemma}

\begin{proof}
Let $\{d_\alpha:\alpha<\mathfrak{d}\}$ be a scale. 
As mentioned in the introduction, this family exists since $\mathfrak{b}=\mathfrak{d}$.

For every $\alpha<\mathfrak{d}$ define $T(\alpha)=T_{d_\alpha}$. Let us see that this family works. Let $\alpha<\beta$. Since $d_\alpha<^*d_\beta$, there exists $k\in\omega$ such that $d_\alpha(m)<d_\beta(m)$ for every $m\geq k$. Let $N\in\omega$ such that $\{s_j:j<k\}\subseteq\bigcup_{i<N}C_i$. Thus $T(\beta)\subseteq T(\alpha)\cup\bigcup_{i<N}C_i$ since $s_m\notin \bigcup_{i<N}C_i$ implies $succ_{T(\beta)}(s_m)\subseteq succ_{T(\alpha)}$.

For (\ref{hechler2}) fix $X\in\textsf{FIN}^{n+1}$ and let $T$ a Hechler tree disjoint from $X$. There is $\alpha<\mathfrak{d}$ such that $f_T<^*d_\alpha$. Pick $k\in\omega$ such that $f_T(m)<d_\alpha(m)$ for every $m\in\omega$ with $s_m\notin\bigcup_{i<k}C_i$. Then
$T(\alpha)\subseteq T\cup\bigcup_{i<k}C_i$ and in consequence $T(\alpha)\cap X\subseteq\bigcup_{i<k}C_i$. Let $\beta>\alpha$. By (\ref{hechler1}), $T(\beta)\subseteq T(\alpha)\cup\bigcup_{i<N'}C_i$. Thus $X\cap T(\beta)\subseteq\bigcup_{i<N}C_i$ for $N=\max\{k,N'\}$.

For the last observation let $H(\alpha)=\omega^{n+1}\setminus T(\alpha)$ and for every $N\in\omega$ let $H(\alpha,N)=H(\alpha)\cup\bigcup_{i<N}C_i$. Therefore $\{H(\alpha,N):\alpha<\mathfrak{d}\land\ N\in\omega\}$ is a cofinal family in $\textsf{FIN}^{n+1}$.
\end{proof}

For the last part of the previous lemma we do not need $\mathfrak{b}=\mathfrak{d}$, in general we have the following:

\begin{corollary}
\emph{$\textsf{cof}(\textsf{FIN}^{n+1})=\mathfrak{d}$.}
\end{corollary}

\begin{proof}
    By the proof of part \ref{hechler2} in the previous Lemma, we have a cofinal family in $\textsf{FIN}^{n+1}$ of size $\d$. 
    Thus $\textsf{cof}(\textsf{FIN}^{n+1})\leq\d$.

    For the other inequality let $\kappa<\d$ and $\{X_\alpha:\alpha\in\kappa\}\subseteq\textsf{FIN}^{n+1}$. For each $X_\alpha$, find a Hechler tree $T_\alpha\subseteq\omega^{n+1}$ such that $X_\alpha\cap T_\alpha=\emptyset$. We can find a function $f\in\w^\w$ that dominates the family $\{f_{T_\alpha}:\alpha<\kappa\}$. Thus $X=\omega^{n+1}\setminus T_f\in\textsf{FIN}^{n+1}$ satisfies that $X\nsubseteq X_\alpha$ for all $\alpha<\kappa$ and the family $\{X_\alpha:\alpha<\kappa\}$ is not cofinal.
\end{proof}

Recall that we identify $\omega^{n+1}$ with the set of increasing sequences of length $n+1$ and given $B\in[\omega]^\omega$, define $B^\uparrow_{n+1}=\{(a_0,\ldots,a_n):\forall i\leq n\ (a_i\in B\ \land\ (a_i<a_{i+1}))\}$. For $i,k\in\omega$ define $C(i,k)=\{(a_0,\ldots,a_n)\in C_i:a_n<k\}$.

\begin{theorem}
\textnormal{$(\diamondsuit(\mathfrak{b})+\mathfrak{d}=\omega_1)$}. For every $n\geq 1$, there exists an $n$-\Ramsey non $(n+1)$-\Ramsey space of size $\w_1$.
\end{theorem}

\begin{proof}
Fix $n\in\omega$ and let $\{T(\alpha):\w\leq\alpha<\omega_1\}$ as in Lemma \ref{hechlertrees}. We will recursively construct an almost disjoint family $\{A_\alpha:\alpha<\omega_1\}\subseteq\omega^{n+1}$. 
For every infinite $\alpha<\omega_1$ let $\alpha=\{\alpha_m:m\in\omega\}$. We will define a Borel function $F:2^{<\omega_1}\to \mathbb{P}$. By a suitable coding we can assume that the domain of $F$ is the set of pairs $(B,\{A_\beta:\beta<\alpha\})$ where
\begin{itemize}
    \item $B\in[\omega]^\omega$,
    \item $\alpha$ is indecomposable,
    \item $A_\alpha\subseteq T(\alpha)$ for $\alpha\geq\omega$,
    \item $\omega\leq\alpha<\omega_1$,
    \item $A_n=C_n$ for every $n\in\omega$ and
    \item $\{A_\beta:\beta<\alpha\}\subseteq[\omega]^\omega$ is an almost disjoint family.
\end{itemize}

Given $(B,\{A_\beta:\beta<\alpha\})$, for every $i\in B\cap succ_{T(\alpha)}(\emptyset)$ define 
$$F(B,\{A_\beta:\beta<\alpha\})(i)=\min\Bigg\{k\in\omega:C(i,k)\cap T(\alpha)\cap B^\uparrow_{n+1}\setminus\left(\bigcup_{m<i}A_{\alpha_m}\right)\Bigg\}.$$

It is not hard to see that $F$ is a Borel function. Let $g:\omega_1\to\omega^\omega$ be a $\diamondsuit(\mathfrak{b})$-sequence with the incarnation of $\mathfrak{b}$ mentioned earlier. 
Now for every $n\in\w$ let $A_n=C_n$ and for every infinite $\alpha<\omega_1$ define 
$$A_\alpha=\bigcup_{i\in\omega}\bigg(C(i,g(\alpha)(i))\setminus\bigcup_{m<i}A_{\alpha_m}\bigg)  \cap T(\alpha).$$
We can assume that every $A_\alpha$ is infinite by increasing $g(\alpha)$ if necessary. It is clear from the definition that $\cA=\{A_\alpha:\alpha<\omega_1\}$ is an almost disjoint family. 
We will now see that $\cF(\cA)$ is $n$-\RRamsey. If $f:[\omega]^n\to\omega^{n+1}$, there exists $X\in[\omega]^\omega$ such that $f''[X]^n\in\textsf{FIN}^{n+1}$ by Lemma \ref{existe un infinito con imagen en fin} and by Lemma \ref{hechlertrees} part \ref{hechler2}, there exists $\alpha<\omega_1$ such that for every $\beta>\alpha$ there is $N\in\omega$ such that $f''[X]^n\cap A_\beta\subseteq f''[X]^n\cap T(\beta)\subseteq\bigcup_{i<N}C_i$. Since $A_\beta\cap C_i\subseteq C(i,g(\beta)(i))$ for every infinite $\beta$, this set is finite and $f''[X]^n$ is almost disjoint with all but countable many $A_\beta$. By Corollary \ref{closure}, $\cF(\cA)$ is $n$-\RRamsey.

To see that it is not $(n+1)$-\Ramsey let $G:[\omega]^{n+1}\to\omega^{n+1}$ where $G$ maps each set $s$ to its increasing enumeration. Thus $G''[B]^{n+1}=B^\uparrow_{n+1}$ for every $B\in[\omega]^\omega$. We will omit the subindex $n+1$. Let $B\in[\omega]^\omega$, we have to show that $B$ is not a convergent subsequence for $G$. This is equivalent to show that there exists $\alpha<\omega_1$ such that $A_\alpha\cap B^\uparrow$ is infinite by Lemma \ref{noramsey}.

Given $(B,\{A_\alpha:\alpha<\omega_1\})$ coded by $f\in2^{\omega_1}$, we can assume that $\{B,\{A_\beta:\beta<\alpha\}\}$ is coded by $f\rest\alpha$ by passing to the club of idecomposable ordinals in $\omega_1$. Since $(B,\{A_\alpha:\alpha<\omega_1\})$ is coded by some $f\in2^{\omega_1}$, there exists $\alpha<\omega_1$ infinite such that $$g(\alpha)\nleq^* F(B,\{A_\beta:\beta<\alpha\}).$$

We use $h$ as a shorthand for $F((B,\{A_\beta:\beta<\alpha\}))$. Hence there are infinitely many $i\in B$ such that $g(\alpha)(i)>h(i)$. Let $B'$ be this set and thus $C(i,h(i))\subseteq C(i,g(\alpha)(i))$ for every $i\in B'$. By the definition of $h$, there exists 
$$b_i\in C(i,h(i))\cap T(\alpha)\cap B^\uparrow_{n+1}\setminus\left(\bigcup_{m<i}A_{\alpha_m}\right)$$
$$\subseteq C(i,g(\alpha)(i))\cap T(\alpha)\setminus\left(\bigcup_{m<i}A_{\alpha_m}\right)\subseteq A_\alpha.$$
As the family $\{C(i,g(\alpha)(i)):i\in B'\}$ is pairwise disjoint, $\{b_i:i\in B'\}$ is infinite and contained in $B^\uparrow\cap A_\alpha$.
\end{proof}

Another property of the parametrized diamonds principles is that they turn the corresponding cardinal invariant small. For instance, $\diamondsuit(\d)\Rightarrow\d=\omega_1$. Thus we can state the previous theorem under $\diamondsuit(\d)$. In Section \ref{section b=s}, we will produce the same examples under $\s=\b$. Although $\diamondsuit(\d)$ (or even $\d=\omega_1$) implies $\s=\b$, a feature of the space constructed in this section is that it has size $\omega_1$ which does not follow from Theorem \ref{teo de s=b}.

\begin{question}
Does $\diamondsuit(\b)$ imply that for every $n\geq 1$ there is an almost disjoint family $\cA$ such that $\cF(\cA)$ is $n$-\Ramsey but not $(n+1)$-\RRamsey?
\end{question}

\section{High dimensional splitting-like cardinal invariants}

We say that $\mathcal{P}=\{P_n:n\in\omega\}\subseteq[\omega]^{<\omega}$ is an {\it interval partition} if $\mathcal{P}$ is a partition of $\omega$ consisting of consecutive intervals.

\begin{mydef}
Let $\mathcal{P}=\{P_n:n\in\omega\}$ be an interval partition and $S\subseteq\omega$. We say that $S$ {\it block-splits} $\mathcal{P}$ if 
$$|\{n\in\omega: P_n\subseteq S\}|=\omega=|\{n\in\omega: P_n\cap S=\emptyset\}|.$$
A family $\mathcal{S}\subseteq[\omega]^\omega$ is a {\it block-splitting} family if for every interval partition $\mathcal{P}$, there exists $S\in\mathcal{S}$ such that $S$ block-splits $\mathcal{P}$.\\
We also define
$$\mathfrak{bs}=\min\{|\mathcal{S}|:\mathcal{S}\textnormal{  is a block splitting family}\}.$$
\end{mydef}

It is easy to see that every block-splitting family is a splitting family. It follows that $\s\leq\b\s$. In fact, we have the following theorem due to Kamburelis and W\k{e}glorz.

\begin{theorem}{\textnormal{\cite{kamburelis1996splittings}}}
$\mathfrak{bs}=\max\{\mathfrak{b,s}\}$
\end{theorem}

Recall that for any $s\in\omega^{\leq n}$, we denote by $C_s$ the set of $t\in\omega^n$ such that $t\upharpoonright dom(s)=s$. If $n\geq 1$ let us denote by $\mathcal{L}_n$ the set of functions $f:[\omega]^n\to\omega^{n+1}$ such that 
there is no $X\in[\omega]^\omega$ so that $f''[X]^n\subseteq C_s$ for some $s\in\omega^{<n+1}\setminus\{\emptyset\}$ (note that it is equivalent if we require this for $s\in\omega^1$).

\begin{mydef}
Let $n\geq 1$, $f\in\mathcal{L}_n$ and $S\in[\omega]^\omega$. We say that $S$ \emph{splits} $f$ (or more precisely {\it $\mathcal{L}_n$-splits} $f$) if there are disjoint $X,Y\in[\omega]^\omega$ such that:
\begin{enumerate}
    \item $f''[X]^n\subseteq S^{n+1}$ and
    \item $f''[Y]^n\subseteq (\omega\setminus S)^{n+1}$
\end{enumerate}

We define $\mu_n$ as the minimum size of an $\mathcal{L}_n$-splitting family, i.e., the minimum size of a family $\mathcal{S}\subseteq[\omega]^\omega$ such that every $f\in\mathcal{L}_n$ is $\mathcal{L}_n$-split by an element in $\mathcal{S}$. 
\end{mydef}

At the moment, it is not even clear that there are $\cL_n$-splitting families, but we will show later in this section that they do exist and find an upper bound for $\mu_n$.

We can consider $\omega^n$ as the set of increasing sequences of $\omega$ of length $n$ through the computable bijection $\varphi$ that maps $(i_0,\ldots,i_{n-1})$ into 
$$\left(i_0\ ,\  i_0+i_1+1\ ,\ \ldots,\sum_{j<n}i_j+(n-1)\right).$$

\begin{notation}
For a set $X\subseteq\omega$ and $n\in\omega$ let 
\begin{itemize}
    \item $nX=\{nx:x\in X\}$ and
    \item $\frac{X}{n}=\{\frac{x}{n}:x\in X\}$ if $n\neq 0$.
\end{itemize}
\end{notation}

We start by finding a lower bound for this new invariants.

\begin{proposition}
$\mathfrak{s}\leq\mu_1$.
\end{proposition}

\begin{proof}
Let $\mathcal{S}$ be an $\mathcal{L}_1$-splitting family of size $\mu_1$. With this, we will define a splitting family of size at most $\mu_1$. For every $X\in[\omega]^\omega$, define $f_X:\omega\to \omega^2$ such that if $X=\{x_n:n\in\omega\}$ is the increasing enumeration of $X$, then $f(n)=(2x_n,2x_n+1)$. Clearly, $f_X\in\mathcal{L}_1$.

There exists $S\in\mathcal{S}$
that $\mathcal{L}_1$-splits $f_X$. Hence there are disjoint $Y,Z\in[\omega]^\omega$ such that $f''Y\subseteq S^2$ and $f''Z\subseteq (\omega\setminus S)^2$. In particular, for every $n\in Y$, we know that $f(n)\in S^2$, which implies that $2x_n\in S$. Similarly, $2x_m\in \omega\setminus S$ for every $m\in Z$.

Thus, we can define $S'=S\cap 2\omega$ and $\widehat{S}=\frac{S'}{2}$ for every $S\in\mathcal{S}$. It is now easy to see that $\{\widehat{S}:S\in\mathcal{S}\}$ is a splitting family of size at most $\mu_1$.
\end{proof}

Given $t\in\omega^n$, we write $t=(t_0,\dots,t_{n-1})$ and denote $im^*(t)=\{t_0,\dots,t_{n-1}\}$. Now we will see that the cardinals $\mu_n$ form a non-decreasing sequence.

\begin{proposition}
Let $n\geq 1$. Every $\cL_{n+1}$-splitting family is also an $\cL_n$-splitting family. Hence $\mu_n\leq\mu_{n+1}$ for every $n\geq 1$.
\end{proposition}

\begin{proof}
Let $\cS$ be an $\cL_{n+1}$-splitting family and let $f\in\cL_n$. Define $f'\in\cL_{n+1}$ such that $f'$ maps the cone $C_s$ into $f(s)\conc k$ where $k$ is the minimum $i$ that can be concatenated to the sequence $f(s)$. Formally,
$$f'(s)=f(s\setminus\max(s))\conc (f(s\setminus\max(s))_{n}+1)$$
for every $s\in[\omega]^{n+1}$. Let $S\in\cS$ that $\cL_{n+1}$-splits $f'$ and let $X,Y\in[\omega]^\omega$ witnessing this. Let $a\in [X]^n$ and pick $x\in X$ such that $x>\max(a)$. Since $f'[X]^{n+1}\subseteq S^{n+2}$, and $im^*(f(a))\subseteq im^*(f'(a\cup\{x\}))$, it is now easy to see that $f(a)\in S^{n+1}$. Therefore $f''[X]^n\subseteq S^{n+1}$. Similarly $f''[Y]^n\subseteq (\omega\setminus S)^{n+1}$.
\end{proof}

We now give an upper bound for $\mu_n$ and implicitly show that $\cL_n$-splitting families exists. In order to achieve this, we need to work with equivalence relations on $[\omega]^n$ and apply Erd\"{o}s-Rado theorem.

\begin{theorem}
For every $n\in\omega$, $\mu_n\leq\mathfrak{bs}$.
\end{theorem}

\begin{proof}
Fix $n\in\omega$. We will show that every block splitting family is $\cL_n$-splitting. Given $f\in\cL_n$ and $a\in[\omega]^n$, we write $f(a)_i=b_i$ if $f(a)=(b_0,b_1,\ldots,b_n)$. Using Corollary \ref{corodeerdosrado}, we find $M\in[\omega]^\omega$ and $I_j\subseteq n$ for every $j\leq n$, such that for every $a,a'\in[M]^n$, $f(a)_j=f(a')_j$ iff $a\rest I_j=a'\rest I_j$.

We claim that $I_j\neq\emptyset$ for every $j\leq n$. Indeed, if $I_j=\emptyset$ for some $j$, there exists $k\in\omega$ such that $f(a)_j=k$ for every $a\in[M]^n$. Since $f(a)$ is increasing for every $a\in[\omega]^n$, there is a finite set $F$ such that $f(a)_0\in F$ for every $a\in[M]^n$, thus applying Ramsey's theorem, we find $M'\in[M]^\omega$ and $l<k$ such that $f(a)_0=l$ for every $a\in[M']^n$. That is, $f''[M']^n\subseteq C_l$, contradicting that $f\in\cL_n$.

For every $j\leq n$, let $m_j=\max(I_j)$. As none of the $I_j$ is empty, the set $\{f(a)_j:a\in[M]^n\}$ is infinite. Moreover,

\begin{itemize}
    \item[($\ast$)]\label{asterisco} for every $j\leq n$ and $k\in\omega$, there exists $z\in M$ such that $f(a)_j>k$ for every $a\in[M]^n$ with $z\in a=\{a_0,\dots,a_{n-1}\}$ and $z=a_i$ for some $i\leq m_j$.
\end{itemize}

To see this, fix $j\leq n$ and $k\in\omega$. For every $k'\leq k$ there exists $a(k')\in[M]^n\cup\{\emptyset\}$ such that $f(a)_j=k'$ if and only if $a\rest I_j=a(k')\rest I_j$ for every $a\in[M]^n$ (possibly, $a(k')=\emptyset$ if $f(a)_j\neq k'$ for every $a\in[M]^n$). Now let $F=\bigcup_{k'\leq k}a(k')$. Thus any $z>\max(F)$ works.

It is also worth noting that given $j\leq n$, $z\in M$ and $t=(t_0,\ldots,t_{m_j-1})\in (M\cap z)^{m_j}$, there is $r(j,z,t)\in\omega$ such that 
\begin{itemize}
    \item[$(\ast\ast)$]\label{dobleasterisco} $f(t\conc \{z\}\conc s)_j=r(j,z,t)$
\end{itemize} 
for every $s\in[M\setminus(z+1)]^{n-m_j}$ (where $t\conc s$ is the concatenation of the sequences $t$ and $s$).

We are going to define an interval partition $\cP=\{P_l:l\in\omega\}$ such that any $S$, that block splits $\cP$, will also $\cL_n$-split $f$. It is clear that if we achieve this, then $\mu_n\leq\mathfrak{bs}$.

To define this partition $\cP$, we will also define an increasing sequence $\{x_l:l\in\omega\}\subseteq M$ and an increasing sequence $\{k_l:l\in\w\}$. 

Let $x_0=\min(M)$ and $k_0=0$. 
For $l>0$, by $(\ast)$, we can find $x_l>\max\{x_i:i<l\}$ in $M$ such that $f(a)_j>k_{l-1}$ for any $j\leq n$ and any $a\in[M]^n$ such that $a_i=x_l$ for some $i\leq m_j$. Given $j\leq n$ and $t=\{t_0,\ldots,t_{m_j-1}\}\subseteq\{x_i:i<l\}$ let $r(j,x_l,t)$ as in $(\ast\ast)$ and define $$k_l>\max(\{r(j,x_l,t):j\leq n\ \land\  t\in[\{x_i:i<l\}]^{m_j}\}\cup\{k_i:i<l\}).$$
Define $P_l=(k_{l-1},k_l]$ for every $l\in\omega$ and let $\cP=\{P_l:l\in\omega\}$.

We proceed to show that $\cP$ is as required. Let $S\in[\omega]^\omega$ that block-splits $\cP$ and define $X=\{x_l:P_l\subseteq S\}$ and $Y=\{x_l:P_l\cap S=\emptyset\}$. Let $a=\{x_{l_i}:i<n\}\subseteq X$ where $x_{l_0}<x_{l_1}<\cdots<x_{l_{n-1}}$. We will prove that $f(a)\in S^{n+1}$. Fix $j\leq n$ and let $i=m_j<n$. By the choice of $l_i$, it follows from the definition of $x_{l_i}$ that $f(a)_j>k_{l_i-1}$. On the other hand, we have that 
$$f(a)_j=f(t\conc x_{l_i}\conc s)_j=r(j,x_{l_i},t)<k_{l_i}$$
with $t=\{x_{l_0},\ldots,x_{l_{m_j-1}}\}$ and $s=\{x_{l_{m_j+1}},\ldots,x_{l_{n-1}}\}$. Thus $f(a)_j\in P_{l_i}$ and $P_{l_i}\subseteq S$ since $x_{l_i}\in X$. As this is true for every $j\leq n$, we have that $f(a)\subseteq S^{n+1}$ and therefore $f''([X]^n)\subseteq S^{n+1}$. A similar argument shows that $f''([Y]^n)\subseteq (\omega\setminus S)^{n+1}$.
\end{proof}

It follows from the results in this section that $\s=\mu_n=\b$ is consistent for every $n\geq 1 $. We do not know if in $\textsf{ZFC}$ one of these equalities hold. 

We ask the following natural questions:
    \begin{question}
    Given $n\geq 1$, is it consistent that $\mu_n<\mu_{n+1}$?
    \end{question} 
    \begin{question}It is true that $\mu_n=\mathfrak{bs}$ for all $n\geq 1$?
    \end{question}
    \begin{question}
    Is it consistent that $\s<\mu_0<\mu_1<\dots<\mu_n<\mu_{n+1}<\dots<\b$?
    \end{question}

\section{$n$-\Ramsey spaces from $\mathfrak{s}=\mathfrak{b}$.}\label{section b=s}

In this section we construct special almost disjoint families following the ideas and techniques developed by Shelah in \cite{shelah} and improved by  Mildenberger, Raghavan and Stepr\={a}ns in \cite{dilip} and \cite{onweaklytight}. For a survey of this method see \cite{hruvsak2007completely}, \cite{guzmancompletely}. We first need to prove that every $\cL_n$-splitting family is \textit{everywhere $\cL_n$-splitting}.
\begin{lemma}
Let $\mathcal{S}$ be a $\cL_n$-splitting family, $A\in [\w]^\w$ and $f:[A]^n\rightarrow \w^{n+1}$ such that there is no $X\in[A]^\w$ so that $f''([X^n])\subseteq C_s$ for some $s\in \w^{<n+1}\setminus\{\emptyset\}$. Then there are $X,Y\in [A]^\w$ and $S\in \mathcal{S}$ such that $f''([X]^n)\subseteq S^{n+1}$ and $f''([Y]^n)\subseteq (\omega\setminus S)^{n+1}$.
\end{lemma}
\begin{proof}
    Let $g:\w\rightarrow A$ a bijection an consider the function $\bar{f}:[\w]^n\rightarrow \w^{n+1}$ given by $\bar{f}(\{n,m\})=f(\{g(n),g(m)\})$. It is easy to see that $\Bar{f}\in \cL_n$, so there is $S\in \mathcal{S}$ and disjoint $X,Y\in [\w]^\w$ such that $\Bar{f}''([X]^n)\subseteq S^{n+1}$ and $\Bar{f}''([Y]^n)\subseteq (\omega\setminus S)^{n+1}$. Taking $\Bar{X}=g''X$ and $\Bar{Y}=g''Y$ we get that $\Bar{f}''[X]^n=f''[\Bar{X}^n]$ and $\Bar{f}''[Y]^n=f''[\Bar{Y}^n]$ so we get the result.
\end{proof}

    Recall that a set $Y\subset\w$ is almost monochromatic (or almost homogeneous) for a partition $f:[\w]^n\rightarrow 2$ if there is a finite set $F\subseteq Y$ such that $f$ is constant on $[Y\setminus F]^n$. The cardinal $\mathfrak{par}_n$ denotes the smallest cardinal $\kappa$ such that there is a family of partitions of $[\w]^n$ of size $\kappa$ such that no infinite set is almost monochromatic for all of them simultaneously.

Note that if we consider partitions into some finite number of pieces $k$, instead of $2$ pieces, we obtain the
same cardinal. Moreover, for all $n\geq 2$ we have that $\mathfrak{par}_n=\mathfrak{par}_2$, in fact:
\begin{theorem}
    (\cite{Blass2010}) For each $n\geq 2$, $\mathfrak{par}_n=\min \{\b,\s\}$.
\end{theorem}

Note that if $\s=\b$ then $\mathfrak{par_2}=\mathfrak{par}_n=\s=\b=\b\s=\mu_n$ for every $n\geq 1$.

  If $S\subseteq\w$ then we will follow the notation $S^0:=S$ and $S^1:=\w\setminus S$. We are now ready to prove the main theorem in this section. It is interesting that the combinatorial characterizations of both $\max\{\b,\s\}$ and $\min\{\b,\s\}$ are used in the proof of this theorem.
 
\begin{theorem}\label{teo de s=b}
    ($\s=\b$) For every $n\geq 1$ there exists an almost disjoint family $\cA$ on $\omega^{n+1}$ of size $\c$ such that $\cF(\cA)$ is $n$-\Ramsey but not $(n+1)$-\RRamsey.
\end{theorem}
\begin{proof}
    Let us fix $\mathcal{S}=\{S_\alpha\mid \alpha\in\b\}$ that is splitting and $\cL_n$-splitting and also let us fix the following enumerations:
    \begin{itemize}
        \item $[\w]^\w=\{B_\alpha\mid \w\leq\alpha<\c\}$.
        \item $\cL_n=\{f_\alpha\mid \w\leq\alpha<\c\}$.
    \end{itemize}
    As usual let $G:[\w]^{n+1}\rightarrow \w^{n+1}$ be the increasing enumeration.

    We want to construct the following families:
    \begin{itemize}
        \item $\{A_\alpha\mid\alpha\in\c\}$,
        \item $\{X_\alpha\mid \w\leq\alpha<\c\}\subseteq[\w]^\w$,
        \item $\{\tau_\alpha\mid \w\leq\alpha<\c\},\{\sigma_\alpha\mid \w\leq\alpha<\c\}\subseteq 2^{<\b}$,
    \end{itemize}
    with the following properties:
    
    \begin{enumerate}
        \item\label{1deconstruccionprincipal} $A_i=C_i$ for all $i\in\w$.
        \item\label{2deconstruccionprincipal} If $\beta<\alpha$ then $\tau_\alpha\nsubseteq \tau_\beta$, $\tau_\alpha\nsubseteq \sigma_\beta$, $\sigma_\alpha\nsubseteq \tau_\beta$ and $\sigma_\alpha\nsubseteq \sigma_\beta$.
        \item\label{3deconstruccionprincipal} If $\beta\in dom(\tau_\alpha)$, then $A_ \alpha\subseteq^* \left(S_\beta^{\tau_\alpha(\beta)}\right)^{n+1}$. 
        \item\label{4deconstruccionprincipal} $f_\alpha''[X_\alpha]^n\in\textsf{FIN}^{n+1}$.
        \item\label{5deconstruccionprincipal} If $\alpha\geq\w$, then $A_\alpha$ is block sequence.
        \item\label{6deconstruccionprincipal} $\forall\alpha\in\c\setminus \w(A_\alpha\subseteq G''[B_\alpha]^{n+1})$.
        \item\label{7deconstruccionprincipal} $\forall \beta,\alpha\in\c\ (\w\leq\beta<\alpha\implies\exists m_{\beta,\alpha}\in\w(f_\beta''[X_\beta\setminus m_{\beta,\alpha}]^n\cap A_\alpha=\emptyset))$.
        \item\label{8deconstruccionprincipal} $\beta\in dom(\sigma_\alpha)\implies\exists n_\beta\in\w(f_\alpha''[X_\alpha\setminus n_\beta]^n\cap \left(S_\beta^{1-\sigma_\alpha(\beta)}\right)^{n+1}=\emptyset)$.
        \item\label{9deconstruccionprincipal} $\forall \beta,\alpha\in \c(\beta\neq\alpha\implies A_\beta\cap A_\alpha=^*\emptyset)$.
    \end{enumerate}

    Let us suppose that we have already constructed the desired objects until step $\alpha$.  We are going to define $X_\alpha$ and $\sigma_\alpha$.\\

\textbf{-$X_\alpha$ and $\sigma_\alpha$ construction:}\\

We are going to construct $\{\alpha_s\mid s\in2^{<\w}\}\subseteq \b$, $\{\eta_s\mid s\in2^{<\w}\}\subseteq 2^{<\b}$, $\{Y_s\mid s\in2^{<\w}\}\subseteq [\w]^\w$, $\{f_s\mid s\in2^{<\w}\}$ such that for all $s\in 2^{<\w}$:
\begin{enumerate}[label=(\alph*)]
    \item\label{Xconstructioncondition1} $f_\emptyset=f_\alpha$, $Y_\emptyset=\w$.
    \item\label{Xconstructioncondition2} 
    $f_s:[Y_s]^n\rightarrow \w^{n+1}$.
    \item\label{Xconstructioncondition3} $\alpha_s=dom(\eta_{s})$. 
    \item\label{Xconstructioncondition4} $Y_{s\conc 0},Y_{s\conc 1}\in[Y_s]^\w$ are disjoint.
    \item\label{Xconstructioncondition5} $f_s''[Y_{s\conc i}]^n\subseteq (S_{\alpha_s}^i)^{n+1}$ for every $i\in 2$.
    \item\label{Xconstructioncondition6} If $\beta\in dom(\eta_{s})$ and $i\in 2$, then $f_{s\conc i}''[Y_{s\conc i}\setminus m]^{n}\cap \left(S_\beta^{1-\eta_s(\beta)}\right)^{n+1}=\emptyset$ for some $m\in\w$.
    \item $s\subseteq t\Rightarrow\eta_s\subseteq\eta_t$.
    \item $s\perp t\Rightarrow\eta_s\perp\eta_t$.
    \item $f_{s\conc i}=f_s\rest Y_{s\conc i}$.
    \item $\alpha_s<\alpha_{s\conc i}$.
\end{enumerate}


    As $\cS$ is $\cL_n$-splitting, then we can find $\alpha_\emptyset=\min\{\beta\in \b\mid S_\beta \text{ splits } f_\emptyset\}$. Then there exist $Z_0,Z_1\in [\w]^\w$ such that $f_\emptyset''[Z_i]^n\subseteq\left(S_{\alpha_\emptyset}^i\right)^{n+1}$ for $i\in 2$. Now  consider the following family of colorings on $[Z_0\cup Z_1]^n$: for every $\beta<\alpha_\emptyset$ let $g_\beta:[Z_0\cup Z_1]^n\rightarrow 3$ as follows:
    \begin{center}
            $g_\beta(\{x_1,\dots,x_n\}) = \left\{
    	       \begin{array}{ll}
                    0      & \mathrm{if\ } f_\emptyset(\{x_1,\dots,x_n\})\in (S_\beta^0)^{n+1} \\
                    1      & \mathrm{if\ } f_\emptyset(\{x_1,\dots,x_n\})\in (S_\beta^1)^{n+1} \\
                    2 & \text{otherwise} \\
    	       \end{array}
	     \right.$
        \end{center}
        Note that for $i\in 2$, $|\{g_\beta\restriction [Z_i]^n\mid \beta<\alpha_\emptyset\}|\leq|\alpha_\emptyset|<\b=\mathfrak{par}_n$, then there exists $Y_i\in [Z_i]^\w$ such that $Y_i$ is almost monochromatic for $g_\beta$ for every $\beta<\alpha_\emptyset$. It is worth to point out that we are defining $\eta_\emptyset$ in terms of $Y_0=Y_{\emptyset\conc0}$ and $Y_1=Y_{\emptyset\conc1}$.

          \begin{claim}\label{Claimetas}
          For every $\beta<\alpha_\emptyset$ there exists $j\in 2$ such that for every $i\in 2$ we have that $f''[Y_{i}\setminus m]^{n}\cap \left(S_\beta^{1-j}\right)^{n+1}=\emptyset$ for some $m\in\w$.
          \end{claim}
           \begin{proof}
        
        Consider $Y_0$ and $Y_1$ and the coloring $g_\beta$. If $Y_i$ is almost $0$-monochromatic for any $i\in2$ then $Y_{1-i}$ can not be almost $1$-monochromatic, since $S_\beta$ does not split $f$. In this case take $j=0$.
        Hence $1-j=1$ and no matter which color $Y_{1-i}$ takes, we know that there exists $n_\beta\in\w$ such that $f_\emptyset([Y_i\setminus n_\beta]^n)\cap (S_\beta^{1})^{n+1}=\emptyset$ for both $i$ and $i-1$. Analogously, if there is an $i<2$ such that $Y_i$ is $1$-monochromatic, $j=1$ works.
        In case both sets are $2$-monochromatic, any $j\in2$ satisfies the claim.
        \end{proof}

        Now define $\eta_\emptyset:\alpha_\emptyset\rightarrow2$ such that for every $\beta<\alpha_\emptyset$ we have that   $\eta_\emptyset(\beta)=j$, where $j$ is given by the claim.       
        For $f_{(0)}$ and $f_{(1)}$ we just take the restriction $f_{(i)}=f\rest[Y_i]^n$.

        Now suppose that we have already constructed $f_s$, $Y_s$, $\alpha_t$ and $\eta_t:\alpha_t\rightarrow 2$ for $s=t\conc j$ where $j\in 2$ and $t\in 2^{<\w}$. We are going to construct $Y_{s\conc 0},Y_{s\conc 1},\alpha_s$ and $\eta_s$.

        As $f_s:[Y_s]^n\rightarrow \w^{n+1}$, there exists a minimum $\alpha_s\in \b$ such that $S_{\alpha_s}$ splits $f_s$. Note that $\alpha_s>\alpha_t$ since no $\beta\leq\alpha_t$ splits $f_s$. 
        Thus there exist $Z_0,Z_1\in [Y_s]^\w$ disjoint with $f_s''[Z_i]^n\subseteq(S_{\alpha_s}^i)^{n+1}$ for every $i\in 2$.
        
        As before, consider the following family of colorings of $[Z_0\cup Z_1]^n$: for every $\beta\in [\alpha_t,\alpha_s)$ let $g_\beta:[Z_0\cup Z_1]^2\rightarrow 3$ as follows:
    \begin{center}
            $g_\beta(\{x_1,\dots,x_n\}) = \left\{
    	       \begin{array}{ll}
                    0      & \mathrm{if\ } f_s(\{x_1,\dots,x_n\})\in (S_\beta^0)^{n+1} \\
                    1      & \mathrm{if\ } f_s(\{x_1,\dots,x_n\})\in (S_\beta^1)^{n+1} \\
                    2 & \text{otherwise} \\
    	       \end{array}
	     \right.$
        \end{center}
        Note that $|\{g_\beta\mid \beta\in [\alpha_t,\alpha_s)\}|\leq|\alpha_s|<\b=\mathfrak{par}_n$, so there exists $Y_{s\conc i}\in [Z_i]^\w$ such that $Y_{s\conc i}$ is almost-homogeneous for every $\beta\in [\alpha_t,\alpha_s)$. At this point we can define $\eta'_s: [\alpha_t,\alpha_s)\rightarrow 2$ with the proper adaptation of Claim \ref{Claimetas}. Now let $\eta_s=\eta_t\cup \eta'_s$. 
        
        From the definition of $\eta'_s$, we know that  for every $\beta\in[\alpha_t,\alpha_s)$, $$f''[Y_{s\conc i}\setminus m_\beta]^n\cap \left(S_\beta^{1-\eta_s(\beta)}\right)^{n+1}=\emptyset$$
        for some $m_\beta\in\w$ and every $i\in2$.

        On the other hand, we know that $Y_{s\conc 0},Y_{s\conc 1}\subseteq Y_s$ and our induction hypothesis implies that if $\beta\in dom(\eta_t)$, then $f''[Y_s\setminus m]^n\cap \left(S_\beta^{1-\eta_t(\beta)}\right)^{n+1}=\emptyset$ for some $m\in\w$, so this last one condition is also true for $Y_{s\conc 0},Y_{s\conc 1}$.
        Finally, let $f_{s\conc i}= f_s\restriction [Y_{s\conc i}]^2$.

        
        


        
    
    This shows that condition \ref{Xconstructioncondition6} holds and all other conditions are clear from the the construction.

    To finish the construction of $X_\alpha$ and $\sigma_\alpha$, let $\eta_g=\bigcup_{n\in\w} \eta_{g\restriction n}$ for every $g\in 2^\w$. 
    Notice that $\eta_g\in2^{<\b}$ as $\b$ has uncountable cofinality (in fact, it is regular, see \cite{Blass2010}). 
    Furthermore, if $f\neq g$, then $\eta_f$ and $\eta_g$ are incompatible nodes of $2^{<\b}$. Since $\alpha<\c$, we can find $g\in 2^\w$ that satisfies that there is no $\beta<\alpha$ such that $\sigma_\beta$ or $\tau_\beta$ extends $\eta_g$. Let $\sigma_\alpha=\eta_g$. Let $X'_\alpha$ be any pseudointersection of $\{Y_{g\restriction n}\mid n\in\w\}$. It follows that $X'_\alpha$ satisfies (\ref{8deconstruccionprincipal}). By Lemma \ref{existe un infinito con imagen en fin}, we can find $X_\alpha\in[X'_\alpha]^\omega$ that also works for (\ref{4deconstruccionprincipal}).\\

    \textbf{-$A_\alpha$ and $\tau_\alpha$ construction.}\\

    We are going to construct a collection $\{B_s\mid s\in 2^{<\w}\}\subseteq [\w]^\w$, $\{\eta_s\mid s\in 2^{<\w}\}\subseteq 2^{<\b}$ and $\{\alpha_s\mid s\in 2^{<\w}\}$ such that:
        \begin{enumerate}[label=(\roman*)]
            \item\label{1 de A y tao} $B_\emptyset=B_\alpha$.
            \item\label{2 de A y tao} $s\subsetneq t\Rightarrow\eta_s\subsetneq\eta_t$.
            \item\label{3 de A y tao} $s\perp t\Rightarrow\eta_s\perp\eta_t$.
            \item\label{4 de A y tao} $\alpha_s=dom(\eta_s)$.
            \item\label{5 de A y tao} $B_{s\conc 0},B_{s\conc 1}\in [B_s]^\w$ are disjoint.
            \item\label{6 de A y tao} If $\beta<\alpha_s$, then $B_s\subseteq^* S_\beta^{\eta_s(\beta)}$.
        \end{enumerate}
Suppose that we have already constructed $\eta_t$ and $\alpha_t$ for every $t\subsetneq s$ and $B_s$. Let $\alpha_s\in\b$ be the minimum such that $|B_s\cap S_{\alpha_s}^0|=\w=|B_s\cap S_{\alpha_s}^1|$, this is possible since $\cS$ is splitting. Notice also that $\alpha_t<\alpha_s$ whenever $t\subsetneq s$ by \ref{6 de A y tao} and \ref{5 de A y tao}. Let $\eta_s:\alpha_s\rightarrow 2$ such that for all $\beta<\alpha_s$, $|B_s\cap S_\beta^{\eta_s(\beta)}|=\w$. Finally, for every $i\in 2$ let $B_{s\conc i}=B_s\cap S_{\alpha_s}^i$. 
Note that for every $\beta\in\alpha_s$ and for $i\in\{0,1\}$, we have that $B_{s\conc i}$ is almost contained in $S_\beta^{\eta_s(\beta)}$. The latter, the fact that $S_{\alpha_s}$ does not split $B_{s\conc i}$ and the fact that $B_s\subseteq B_t$ for $t\subseteq s$, imply that $\eta_{s}\subsetneq\eta_{s\conc i}$ for every $i\in 2$. 

    Now, for every $g\in 2^\omega$ let $\eta_g=\bigcup_{n\in\w}\eta_{g\restriction n}$. As before $\eta_g\in 2^{<\b}$ and if $f\neq g$, $\eta_f$ and $\eta_g$ are incompatible nodes of $2^{<\b}$. As $\alpha<\c$, we can find $g\in 2^\w$ such that there is no $\beta<\alpha$ such that $\sigma_\beta$ or $\tau_\beta$ extends $\eta_g$, then let $\tau_\alpha=\eta_g$. Consider $Z$ a pseudointersection of $\{B_{g\restriction n}\mid n\in\w\}$ that is contained in $B_\emptyset=B_\alpha$. 
    
    Note that if $\beta<\alpha$ and $\sigma_\alpha$ and $\tau_\beta$ are incompatibles, then there exists $m\in\w$ such that if $\xi=\sigma_\alpha\vartriangle\tau_\beta$\footnote{Recall that $\sigma\vartriangle\tau=\min\{\alpha\mid\sigma(\alpha)\neq\tau(\alpha)\}$} then $f_\beta''[X_\beta\setminus m]^n\cap\left(S_\xi^{\sigma_\alpha(\xi)}\right)^{n+1}=\emptyset$. In particular, as $Z$ is almost contained in $S_\xi^{\sigma_\alpha(\xi)}$, then $f_\beta''[X_\beta\setminus m]^n\cap (Z)^{n+1}=^*\emptyset$.

Consider the following families of elements of $\mathsf{FIN}^{n+1}$:
$$\cH_0=\{f_\beta''[X_\beta]^n\mid \omega\leq\beta<\alpha\land(\tau_\alpha\supseteq \sigma_\beta)\},$$
$$\cH_1=\{A_\beta\mid \omega\leq\beta<\alpha(\tau_\alpha\supseteq \tau_\beta)\}$$
and
$$\cH_2=\{A_i\mid i\in\w\}.$$
Define $\cH=\cH_0\cup\cH_1\cup\cH_2$.
Note that $|\cH|<\b$, 
so we can apply Lemma \ref{cov igual a b} to the family $\cH$ and the infinite set $Z$ to obtain $A_\alpha$ a block sequence of elements of $Z$ such that its intersection with every element of $\cH$ is finite. 

Given $\beta\in dom(\tau_\alpha)$, there is $s\in2^{<\w}$ such that $\beta\in dom(\eta_s)$ and $s\subseteq g$. Thus $Z\subseteq^*B_s\subseteq^*S_\beta^{\eta_s(\beta)}=S_\beta^{\tau_\alpha(\beta)}$. Therefore $A_\alpha\subseteq^*\left(S_\beta^{\tau_\alpha(\beta)}\right)^{n+1}$ as $A_\alpha\subseteq Z^{n+1}$ is a block sequence.

Let us verify that this construction satisfies conditions (\ref{1deconstruccionprincipal}) to (\ref{9deconstruccionprincipal}).
Conditions (\ref{1deconstruccionprincipal}), (\ref{2deconstruccionprincipal}), (\ref{4deconstruccionprincipal}), (\ref{5deconstruccionprincipal}) are true by definition. Condition (\ref{6deconstruccionprincipal}) holds since $A_\alpha$ is a block sequence of elements of $Z\subseteq B_\alpha$. Condition (\ref{8deconstruccionprincipal}) follows from the choice of $\sigma_\alpha$. We have already checked condition (\ref{3deconstruccionprincipal}) at the end of the construction of $A_\alpha$. It remains to show that (\ref{7deconstruccionprincipal}) and (\ref{9deconstruccionprincipal}) hold.

\underline{Condition (\ref{7deconstruccionprincipal}):} Let $\beta$ and $\alpha$ in $\c$ such that $\w\leq \beta< \alpha$. If $\sigma_\beta\subseteq \tau_\alpha$, then for the construction of $A_\alpha$ we know that $f_\beta''[X_\beta]^n\cap A_\alpha=^*\emptyset$ and, as $A_\alpha$ is block sequence, then there exists $m_{\beta,\alpha}\in\w$ such that $f_\beta''[X_\beta\setminus m_{\beta,\alpha}]^n\cap A_\alpha=\emptyset$.

On the other hand if $\sigma_\beta$ and $\tau_\alpha$ are incompatibles then exists $\xi\in dom(\tau_\alpha)\cap dom(\sigma_\beta)$ such that $\tau_\alpha(\xi)=i$ and $\sigma_\beta(\xi)=1-i$.  By (\ref{3deconstruccionprincipal}) and (\ref{8deconstruccionprincipal}) applied to $\xi$ we get the result.


\underline{Condition (\ref{9deconstruccionprincipal}):} It is clear that if $i\neq j$ are finite ordinals then $A_i\cap A_j=\emptyset$. Also, if $\w\leq\alpha<\c$ then $A_\alpha$ is block sequence, so $A_\alpha\cap A_i$ is finite for every $i\in\w$. Now let $\w\leq\beta<\alpha<\c$. If $\tau_\alpha\supseteq \tau_\beta$ then $A_\alpha\cap A_\beta=^*\emptyset$ by the construction of $A_\alpha$, since $A_\beta\in\cH$. Otherwise $\tau_\alpha$ and $\tau_\beta$ are incompatible nodes in $2^{<\b}$, so there exists $\xi\in dom(\tau_\alpha)\cap dom(\tau_\beta)$ such that $A_\alpha\subseteq^*(S_\xi^{\tau_\alpha(\xi)})^{n+1}$ and $A_\beta\subseteq^*(S_\xi^{1-\tau_\alpha(\xi)})^{n+1}$, in particular $A_\alpha\cap A_\beta$ is finite.
This finishes the construction.

Let us see that $\cF(\cA)$ is the desired space. Conditions (\ref{3deconstruccionprincipal}) and (\ref{6deconstruccionprincipal}) imply that for every infinite $\alpha\in\c$, the elements of the almost disjoint family in the closure of $f_\alpha''[X_\alpha]^n$, are at most those $A_\beta$ such that $\tau_\beta\subseteq\sigma_\alpha$. In particular, the closure of $f_\alpha''[X_\alpha]^n$ has size less than $\b$.
By Corollary \ref{lema de paul} we get that $\cF(\cA)$ is $n$-\RRamsey. On the other hand, condition (\ref{7deconstruccionprincipal}) imply that condition \ref{lemanorramsey4} of Lemma \ref{noramsey} holds and the other three conditions when $f=G$ are clear, so $\cF(\cA)$ is not $(n+1)$-\RRamsey.
\end{proof}

Most modifications of the technique developed by Shelah used previously in the literature use assumptions of the form $\iota\leq\theta$ where $\iota$ is a kind of splitting cardinal invariant and $\theta$ is another cardinal invariant imposed by the particularities of the problem. We do not know if the previous result can be weakened to $\s\leq\b$.

\begin{question}
Does it follow from $\s\leq\b$ that there are $n$-\Ramsey spaces that fail to be $(n+1)$-\RRamsey?
\end{question}

Or even better
\begin{question}
Is there in $\mathsf{ZFC}$, a space that is $n$-\Ramsey but no $(n+1)$-\Ramsey for every $n\geq 2$? 
\end{question}

\section{On the character of $n$-\Ramsey spaces}

One of the key facts in the constructions of our examples has been that any sequentially compact space of character less than $\b$ is $n$-\Ramsey for every $n\geq 1$. It is then natural to ask what is the smallest possible character of an $n$-\Ramsey space that is not $(n+1)$-\RRamsey: 

\begin{mydef}
    If $n\geq 2$ let $\lambda_{n}$ be the minimum character of a space $X$ such that $X$ is $(n-1)$-\Ramsey but not $n$-\RRamsey.
\end{mydef}

We do not know if these cardinals are well defined since we do not know if the are $n$-\Ramsey spaces that are not $(n+1)$-\Ramsey in \textsf{ZFC} (for $n>1$). To avoid unnecessary complications, let us agree that $\lambda_n=\c$ when it is not well defined by the previous definition. Theorem \ref{teoremadepaulykubis} gives us $\b$ as a lower bound for $\lambda_n$ for every $n\geq 2$. Now, we will see that, using the following space, introduced by Burke and van Douwen in \cite{burkevandouwen}, we can actually compute $\lambda_2$.

Let $\{f_\alpha\mid\alpha\in\b\}$ consisting of strictly increasing functions, which is well ordered by $\leq^*$ and for each $\alpha\in\b$ let $B_\alpha=\{(n,m)\mid m\leq f_\alpha(n)\}$. Now consider $\cB:=\{B_\alpha\mid\alpha\in\b\}$ and 
let $X=(\w\times\w)\cup \cB$ be a topological space, where the elements of $\w\times\w$ are isolated and a basic neighborhood for $B_\beta\in\cB$ has the form
\[N(B_\alpha,B_\beta,F)=(B_\alpha,B_\beta]\cup ((B_\beta\setminus B_\alpha)\setminus F),\] where $\alpha\in \beta\cup \{-1\}$, $B_{-1}=\emptyset$ and $F\in [\w\times\w]^{<\w}$. It is easy to check that $X$ is Hausdorff, separable and since the basic neighborhoods of every point are compact, then $X$ is also locally compact (see Example 7.3 of \cite{van1984integers}). On the other hand, $X$ is not sequentially compact since for all $n\in\w$, the sequence $((n,m))_{m\in\w}$ does not converge in $X$. So let $\overline{X}=X\cup\{\infty\}$ be the one-point compactification of $X$. Then the basic neighborhoods of $\infty$ have the form:
\[
\{\infty\}\cup(X\setminus C)
\]
where $C\subseteq X$ is compact.

\begin{theorem}\label{lashipotesisdelteodekubisypaulsonoptimas}
        The space $\overline{X}$ is a compact and sequentially compact space of character $\b$ that is not $2$-\RRamsey.
\end{theorem}
\begin{proof}
    Let us see that $\overline{X}$ is sequentially compact. Let $f:\w\rightarrow \overline{X}$. If $f$ takes the value $\infty$ infinitely many times, then we are done. Also if $f$ takes values in $\cB$ infinitely many times, then we are done since $\b$ has uncountable cofinality (and then it is sequentially compact). So let us suppose that $f:\w\rightarrow \w\times\w$. If $f''\w\in\cA$ then there exists $B_\beta\in \cB$ with $\beta\in\b$ minimum such that $f''\w\cap B_\beta$ is infinite. Taking $B\in [\w]^\w$ such that $f''B\subseteq B_\beta$ we have that $f\restriction B$ converges to $B_\beta$. If $f''\w\not\in \cA$ then let $B\in [\w]^\w$ and $k\in\w$ such that $f''B\subseteq A_k$. It follows that $f\restriction B$ converges to $\infty$ since every basic neighborhood of $\infty$ contains all but finite elements of $A_n$ for every $n\in\w$.
    
    We will now see that $\overline{X}$ is not $2$-\RRamsey, so consider $G:[\w]^2\rightarrow\overline{X}$ given by $G(\{n,m\})=(n,m)$ where $n<m$. We will see that for all $B\in [\w]^\w$, $G\restriction [B]^2$ does not converge in $\overline{X}$. So let us fix $B\in [\w]^\w$ and $B_\beta\in\cB$. Note that for every basic neighborhood $N$ of $B_\beta$ we have that $N\cap A_n$ is finite for every $n\in\w$, but $G''[B\setminus F]^2$ has infinite intersection with every $A_n$ such that $n\in B\setminus F$. This shows that there is no $F\in [\w]^{<\w}$ such that $G''[B\setminus F]^2\subseteq N$. On the other hand, as $G$ is injective, $G\restriction [B]^2$ does not converge to any isolated point. So the only thing that remains to be seen is that $G\restriction [B]^2$ does not converge to $\infty$. 

    Let $(m_k)_{k\in\w}$ be the increasing enumeration of $B$ and let $A=\{(m_k,m_{k+1})\mid k\in\w\}$. Note that:
    \begin{enumerate}\label{estrella}
    \item\label{1deestrella} $A\in \cA$.\item\label{2deestrella} For all $F\in[\w]^{<\w}$ we have that $A\subseteq^* G''[B\setminus F]^2$. 
    \end{enumerate}
     Now let $B_\beta\in\cB$ such that $|B_\beta \cap A|=\w$ and $\beta$ is the minimum with this property. Consider the following basic neighborhood of $\infty$:
    \[
    N=\{\infty\}\cup (X\setminus N(B_{-1},B_\beta,\emptyset))=\{\infty\}\cup \left(\cB\setminus (B_{-1},B_\beta ]\right)\cup \left((\w\times\w) \setminus B_\beta\right). 
    \]
    Remember that $B_{-1}=\emptyset$. If there exists $F\in [\w]^{<\w}$ such that $G''[B\setminus F]^2\subseteq N$, then in particular we have that $G''[B\setminus F]^2\cap B_\beta=\emptyset$, and as $A\subseteq^* G''[B\setminus F]^2$, then $A\cap B_\beta$ is finite, which is a contradiction. This implies that there is no such $F$ and then $G\restriction [B]^2$ does not converge, which proves that $\overline{X}$ is not $2$-\RRamsey.  
    \end{proof}

With this last result and Theorem \ref{teoremadepaulykubis}, we get the following nice topological characterization of $\b$:
\begin{corollary}
$\b=\lambda_2$, i.e.:
    $$\b=\min\{\chi(X)\mid X \text{ is sequentially compact but not $2$-\RRamsey}\}.$$
\end{corollary}

We can still ask for the minimal size of an almost disjoint family $\cA$ such that $\cF(\cA)$ is not $2$-\RRamsey, i.e., the minimal character of a space of the form $\cF(\cA)$ that is not $2$-\RRamsey. In general we can ask this for every $n\geq2$:

\begin{mydef}
    If $n\geq 2$ let \emph{$\a_{n}^{sc}$} be the minimum size of an almost disjoint family such that $\cF(\cA)$ is $(n-1)$-\Ramsey but not $n$-\RRamsey.
\end{mydef}

It follows from the results in \cite{kubis2021topological} that $\b\leq\a_{2}^{sc}\leq\a$. Also, it should be clear that $\lambda_n\leq\a_n^{sc}$ for all $n\geq 2$ whenever they are well defined.

\begin{mydef}
\cite{shelah} Let $\a_*$ be the minimal size of an ad family $\cA$, such that there are (disjoint) $C_0, C_1, C_2,\ldots\in\cA^\perp$ so that for every $D\in [\w]^\w$, if $|C_n \cap D|=\w$ for infinitely many $n\in\w$, then there is $A\in \cA$ for which $|A\cap D|=\w$.
\end{mydef}

The following notions were introduced by Arhangel’skii \cite{arhangel1981frequency}:
\begin{mydef}
    Let $X$ be a topological space and $x\in X$. Then $x$ is an \emph{$\alpha_i$-point} ($i=1,2,3$) if for every family $\{S_n\mid n\in\w\}$ of sequences from $X$ converging to $x$, there exists a sequence $S$ that converges to $x$ such that:
    \begin{itemize}
        \item ($\alpha_1$) $S_n\setminus S$ is finite for all $n\in\omega$, 
        \item ($\alpha_2$) $S\cap S_n\neq\emptyset$ for all $n\in\w$.
        \item ($\alpha_3$) $|S\cap S_n|=\w$ for infinitely many $n\in\w$.
    \end{itemize}
    A space X is an $\alpha_i$-space if every point $x\in X$ is an $\alpha_i$-point.
\end{mydef}

Arhangel'skii introduced the $\alpha_i$ properties for $i\in\{1,2,3,4,5\}$ in \cite{arhangel1981frequency}. It is easy to see that all $\alpha_1$-points are $\alpha_2$-points and that all $\alpha_2$-points are $\alpha_3$-points. We refer the reader to \cite{arhangel1981frequency} for the rest of the definitions and a deep study on their relationship.

\begin{proposition}\label{a estrellita y alpha 3}
$\a_*$ is the minimum size of an ad family such that $\cF(\cA)$ is not $\alpha_3$.    
\end{proposition}
\begin{proof} 
    Let $\kappa<\a_*$ and be $\cA$ ad family of cardinality $\kappa$. As we want to prove that $\cF(\cA)$ is $\alpha_3$, it is enough to prove that $\infty$ is an $\alpha_3$ point since any other point has countable character. Let $\{S_n\mid n\in\w\}$ be a collection of convergent sequences to $\infty$. We can assume that $S_n$ only takes values in $\w$ for every $n\in\w$, since any sequence in $\cA$ converges to $\infty$. As $S_n$ is convergent to $\infty$, $S_n\in \cA^\perp$ for every $n\in\w$. Since $|\cA|<\a_*$, there exists $B\in[\w]^\w$ such that $S_n\cap B$ is infinite for infinitely many $S_n$ and $A\cap B$ is finite for every $A\in\cA$. Then $B$ is the required convergent sequence to $\infty$. 

    Conversely, let $\cA$ be an ad family such that $\cF(\cA)$ is $\alpha_3$. We want to prove that $\cA$ is not a witness of $\a_*$, so let $C_0, C_1, C_2,\ldots\in\cA^\perp$. Note that every $C_n$ is a convergent sequence to $\infty$ and as $\infty$ is an $\alpha_3$-point, there exists a sequence $B\subseteq\Psi(\cA)$ that converges to $\infty$ and such that $|C_n\cap B|=\w$ for infinitely many $n\in\w$. 
    We can assume that $B\subseteq\omega$ since otherwise $B'=B\cap(\bigcup_{i\in\omega})C_i$ is also a witness for the $\alpha_3$-property for the family $\{C_i\mid i\in\omega\}$.
    As $B$ converges to $\infty$, we also have that $B\cap A=^*\emptyset$ for all $A\in\cA$, so $\cA$ is not a witness for $\a_*$ as we wanted to show.
\end{proof}
Following the idea of the previous proof, we can actually show that $2$-\Ramsey ad families are $\alpha_3$. Given a non $\alpha_3$ ad family $\cA$ and a family $\{S_n:n\in\omega\}\subseteq\cA^\perp$ that witnesses this, we can define $f:[\omega]^2\to\omega$ by $f(\{n,m\})=S_n(m)$, where $S_n(m)$ is the $m^{\textnormal{th}}$ element of $S_n$. Thus $\cA$ and
$\cA_0=\{S_n\mid n\in\w\}$ satisfy the hypothesis of Lemma \ref{noramsey} for $n=2$ and $\cF(\cA)$ is not $2$-\RRamsey.

In general, it is possible to prove that $2$-\Ramsey spaces are $\alpha_3$. This result was independently proved by Szeptycki, Almontashery and Memarpanahi. Moreover, they also proved that $2$-\Ramsey spaces do not need to be $\alpha_2$ and that $\alpha_1$ spaces need not to be $2$-\RRamsey. 



\begin{theorem}
    If $X$ is $2$-\RRamsey, then $X$ is $\alpha_3$.
\end{theorem}

\begin{proof}
    Let $x\in X$ and for every $n\in\w$ let $Y_n=\{x(n,m):m\in\omega\}$ be a sequence converging to $x$. Define $f:[\omega]^2\to X$ by $f(n,m)=x(n,m)$, 
    where $n<m$. Since $X$ is $2$-\RRamsey, there exists $B\in [\w]^\w$ such that $f\restriction [B]^2$ converges.

    Let $Y:=f''[B]^2\in [X]^\w$. It is clear that $Y\cap f_n''\w$ is infinite for every $n\in B$, so we are done if we prove that $Y$ is convergent to $x$. Let $U$ be an open set such that $x\in U$. Since $f\restriction[B]^2$ converges to $x$, we know that there exists $m\in\w$ such that $f''[B\setminus m]^2\subseteq U$. It follows that $Y\setminus U\subseteq \bigcup_{n\in B\cap m} Y_n$. But $\bigcup_{n\in B\cap m} Y_n$ is a convergent sequence to $x$, so $\bigcup_{n\in B\cap m} Y_n\setminus U$ is finite, then
    \[
    Y\setminus U\subseteq (\bigcup_{n\in B\cap m} Y_n)\setminus U
    \]
    and the last is finite. We conclude that $Y$ converges to $x$.
\end{proof}

It is also worth to mention that the previous proof follows the same idea as the proof that the spaces with the Ramsey property considered by Knaust \cite{angelicspaces} are $\alpha_3$. Recall that a space $X$ has the Ramsey property if whenever $\lim_{i\to\infty}\lim_{j\to\infty}x(i,j)=x$, there exists $B\in[\w]^w$ such that $f\rest [B]^2$ converges to $x$, where $f(\{i,j\})=x(i,j)$.

\begin{corollary}\label{si A es ad 2-ramsey entonces es alpha_3}
    If $\cA$ is an ad family such that $\cF(\cA)$ is $2$-\RRamsey, then $\cF(\cA)$ is $\alpha_3$.
\end{corollary}

In \cite{genericexistenceofmadfamilies} it is proved that $\b\leq\a_*\leq\a$ and Corollary \ref{si A es ad 2-ramsey entonces es alpha_3} implies that $\a_2^{sc}\leq\a_*$, so we have that $\b=\lambda_2\leq\a_2^{sc}\leq \a_*\leq\a$. Of course, the inequality $\lambda_n\leq\a_n^{sc}$ is clear from the definition. This is pretty much what we know so far about these cardinals. 

\begin{question}
    Are there any other relations between the cardinals $\lambda_n$ and $\a_n^{sc}$ with $\b$ and $\a$?
\end{question}

\begin{question}
    Is it consistent that $\a_2^{sc}<\a_*$? 
\end{question}

\begin{question}
    Are all the $\lambda_n$ and $\a_m^{sc}$ consistently different? In particular, is $\a_2^{sc}<\cdots<\a_n^{sc}<\cdots<\a$ consistent?
\end{question}

\subsection*{Acknowledgments}
We would like to thank Paul Szeptycki and Stevo Todor\v{c}evi\'c for all their comments regarding the content of the paper. Part of the research was done while the authors were part of the Thematic Program on Set Theoretic Methods in Algebra, Dynamics and Geometry at the Fields Institute. The first author also wishes to thank York University and the Fields Institute for his Postdoctoral Fellowship.

\bibliography{ref}{}
\bibliographystyle{plain}

\end{document}